\documentclass[11pt,reqno]{amsart}
\usepackage[T1]{fontenc}
\usepackage{lmodern}
\usepackage{amsmath,amssymb,mathtools}
\usepackage[margin=1.1in]{geometry}
\usepackage{microtype}
\usepackage{needspace}
\usepackage[colorlinks=true,linkcolor=blue,citecolor=blue,urlcolor=blue]{hyperref}
\newtheorem{theorem}{Theorem}[section]
\newtheorem{proposition}[theorem]{Proposition}
\newtheorem{lemma}[theorem]{Lemma}
\newtheorem{corollary}[theorem]{Corollary}
\theoremstyle{definition}
\newtheorem{definition}[theorem]{Definition}
\newtheorem*{definition*}{Definition}
\newtheorem{example}[theorem]{Example}
\theoremstyle{remark}
\newtheorem{remark}[theorem]{Remark}
\newtheorem*{remark*}{Remark}
\newcommand{\C}{\mathbb C}

\newcommand{\B}{\mathcal B}
\newcommand{\A}{\mathcal A}
\newcommand{\M}{\mathcal M}

\newcommand{\F}{\mathcal F}
\newcommand{\bF}{\mathbb F}
\makeatletter
\def\moverlay{\mathpalette\mov@rlay}
\def\mov@rlay#1#2{\leavevmode\vtop{%
\baselineskip\z@skip \lineskiplimit-\maxdimen
\ialign{\hfil$#1##$\hfil\cr#2\crcr}}}
\makeatother
\def\plangle{\moverlay{(\cr<}}
\def\prangle{\moverlay{)\cr>}}
\newcommand{\freefield}[1]{\C\mathopen{\plangle}#1\mathclose{\prangle}}
\newcommand{\Div}{\operatorname{Div}}
\newcommand{\rank}{\operatorname{rank}}
\newcommand{\dom}{\operatorname{Dom}}

\newcommand{\Image}{\operatorname{Im}}

\newcommand{\dist}{\operatorname{dist}}
\title[Rationality criteria]{Rationality criteria of affiliated operators}
\author{Sheng Yin}
\address{Institute for Advanced Study in Mathematics, Harbin Institute of Technology}
\email{sheng.yin@hit.edu.cn}
\date{\today}
\begin{document}
\begin{abstract}
Let $X:=(X_1,\ldots,X_d)$ be a tuple in a tracial $W^*$-probability space.
We establish three criteria for affiliated operators to be rational in $X$, each with several equivalent characterizations.
The first criterion is of algebraic nature and assumes only $\Delta(X)=d$ for $X$.
From it, we derive a second criterion of an analytic nature using liberation dual systems.
When $X$ consists of normal generators and admits a dual system, we obtain a third criterion that further refines the second in terms of finite-rank fraction commutators, thereby extending the rationality criteria of Duchamp and Reutenauer, Linnell, and Miyagawa for affiliated operators.
\end{abstract}
\maketitle
\section{Introduction}\label{sec:introduction}

Evaluating a noncommutative rational function at a tuple $X$ of bounded operators can yield an unbounded affiliated operator.
The free-field realization theorem of Mai, Speicher and Yin \cite[Theorem 1.1]{MSY23} shows that evaluation at $X$ embeds the free field into the $*$-algebra of affiliated operators when $\Delta(X)=d$.
Under this condition, we first establish algebraic characterizations of rationality in $X$.
We then derive two analytic criteria: one using liberation dual systems and the other using dual systems.

Our question is motivated by Connes's rationality conjecture for free group algebras \cite[Section IV.5, Remark 3]{Con94}, a noncommutative analogue of Kronecker's criterion \cite{Kro81}.
In its operator formulation, this criterion characterizes rationality by a finite-rank commutator condition \cite[Section IV.3, Proposition 1]{Con94}.
For a free group, Connes observed that commutators with the operator arising from its Cayley tree have finite rank on the rational closure of the group algebra in the reduced group $C^*$-algebra,
and conjectured the converse.

Duchamp and Reutenauer \cite[Theorem 7]{DR97} proved this conjecture for a free group $\bF_d$ with $d$ generators ($d\geq1$),
and showed that the rational closure coincides with the division closure.
Their proof uses noncommutative rational series.
Using the bounded dual system for free Haar unitaries (see, for example, \cite[Proposition 6.3]{MSY23}),
one can express the Cayley-tree finite-rank condition as finite rank of the corresponding commutators.

Linnell \cite[Theorem 1.1]{Lin00} extends the free-group criterion to affiliated operators using left and right fraction representations.
Theorem~\ref{thm:analytic} extends this criterion to normal generators admitting a dual system in the sense of Definition~\ref{def:dual-system}, including free Haar unitaries with their bounded dual systems.
Linnell also identifies the bounded part of the affiliated division closure with the division closure in the reduced group $C^*$-algebra.
Together with this identity, Theorem~\ref{thm:analytic} recovers the criterion of Duchamp and Reutenauer.

Miyagawa \cite[Theorems 3.1 and 4.5]{Miy22} proves a semicircular analogue of these criteria, characterizing rationality by finite-rank commutators with the right annihilation operators (whose negatives form a bounded dual system in our convention).
His proof adapts the noncommutative rational-series approach of Duchamp and Reutenauer to the semicircular setting.
He also extends this characterization to affiliated operators using left and right fraction representations.

\subsection{A rationality criterion under maximal \texorpdfstring{$\Delta$}{Delta}}
Our first main result gives several equivalent characterizations under maximal $\Delta$ alone.
The quantity $\Delta$ was introduced by Connes and Shlyakhtenko in their paper \cite{CS05} on the homology of von Neumann algebras.

Throughout the paper, $(\M,\tau)$ is a tracial $W^*$-probability space, that is, a von Neumann algebra $\M$ equipped with a faithful normal tracial state $\tau$.
Let $X=(X_1,\ldots,X_d)\in\M^d$, with $d\geq1$.
Put $H:=L^2(\M,\tau)$ and $\Omega:=1$.
Inner products on $H$ are linear in the second variable, and we identify elements of $\M$ with left multiplication on $H$.
Write $\F(H)$ for the space of finite-rank operators on $H$.
For $x\in\M$, let $\rho(x)$ denote the right multiplication operator on $H$ given by $\rho(x)(y\Omega):=(yx)\Omega$ for $y\in\M$.
The condition $\Delta(X)=d$ means that, for $T_1,\ldots,T_d\in\F(H)$,
\[
 \sum_{i=1}^d[\rho(X_i),T_i]=0
 \quad\Longrightarrow\quad T_1=\cdots=T_d=0;
\]
see Subsection~\ref{subsec:delta-tensor-maps} for details.

We write $L^0(\M)$ for the $*$-algebra of closed densely defined operators affiliated with $\M$.
We identify $H$ with the square-integrable affiliated operators in $L^0(\M)$ \cite[Theorem~7.3.2]{AP18}.
In particular, $H\otimes_{\C}H\subseteq L^0(\M)\otimes_{\C}L^0(\M)^{\mathrm{op}}$.
Write $\A_X:=\C\langle X_1,\ldots,X_d\rangle$ for the algebra generated by $X_1,\ldots,X_d$.
We denote its division closure in $L^0(\M)$ by $\Div_{L^0(\M)}\A_X$ and call an affiliated operator \emph{rational} when it belongs to this closure; see Definition~\ref{def:rationality}.

\begin{definition}\label{def:fraction-pair}
For $u\in L^0(\M)$, a \emph{fraction pair} is a quadruple $Q=(s,a,b,t)\in\M^4$ satisfying
\begin{equation}
 u=s^{-1}a=bt^{-1},\qquad s,t\text{ invertible in }L^0(\M).
 \label{eq:fractions}
\end{equation}
\end{definition}

\begin{definition}\label{def:fraction-commutator}
Let $D$ be a linear operator on $H$ with domain $\dom D$.
For a fraction pair $Q=(s,a,b,t)$ of an element $u\in L^0(\M)$, define the \emph{fraction commutator of $D$ and $Q$} by
\begin{equation}
 \Gamma[Q,D]:=sDb-aDt,
 \label{eq:fraction-commutator}
\end{equation}
with domain given by $\{\xi\in H: b\xi,t\xi\in\dom D\}$.
\end{definition}
The identities $ut=b$ and $su=a$ from \eqref{eq:fractions} give formally
\[
 s[D,u]t=sD(ut)-(su)D t=sD b-aD t=\Gamma[Q,D].
\]
That is, $s$ and $t$ remove the fractions in $u$ from the commutator expression.

Multiplication in $L^0(\M)\otimes_{\C}L^0(\M)^{\mathrm{op}}$ is given by
\[
 (u_1\otimes v_1^{\mathrm{op}})(u_2\otimes v_2^{\mathrm{op}})
 :=u_1u_2\otimes(v_2v_1)^{\mathrm{op}}.
\]
To state our first main result, we use the right multiplication maps in \cite[Section~3.1.1, equation~(3.1)]{CS05},
with the coefficient algebra enlarged from $\M\otimes_{\C}\M^{\mathrm{op}}$ to $L^0(\M)\otimes_{\C}L^0(\M)^{\mathrm{op}}$.

\begin{definition}\label{def:psi}
For $w\in L^0(\M)$, define the right multiplication map
\[
 \Psi_w:L^0(\M)\otimes_{\C}L^0(\M)^{\mathrm{op}}
 \longrightarrow L^0(\M)\otimes_{\C}L^0(\M)^{\mathrm{op}}
\]
by
\begin{equation}
 \Psi_w(k):=k(w\otimes1-1\otimes w^{\mathrm{op}})
 \label{eq:psi}
\end{equation}
\end{definition}

For $u,v\in L^0(\M)$, this gives
\[
 \Psi_w(u\otimes v^{\mathrm{op}})
 =uw\otimes v^{\mathrm{op}}-u\otimes(wv)^{\mathrm{op}}.
\]
For a tuple $X=(X_1,\ldots,X_d)\in\M^d$, we further write $\Psi_i:=\Psi_{X_i}$.

\Needspace{4\baselineskip}
\begin{theorem}\label{thm:delta}
Suppose $\Delta(X)=d$.
For $u\in L^0(\M)$ the following are equivalent:
\begin{enumerate}
\item $u\in\Div_{L^0(\M)}\A_X$.
\item\label{item:delta-tensors} There are $k_i\in L^0(\M)\otimes_{\C}L^0(\M)^{\mathrm{op}}$ such that
\begin{equation}
 \sum_{i=1}^d\Psi_i(k_i)=u\otimes1-1\otimes u^{\mathrm{op}}.
 \label{eq:tensor-criterion}
\end{equation}
\item\label{item:delta-all-fraction-tensors}
For every fraction pair $Q=(s,a,b,t)$ for $u$, there exist $\ell_1,\ldots,\ell_d\in L^0(\M)\otimes_{\C}L^0(\M)^{\mathrm{op}}$ such that
\begin{equation}\label{eq:all-fraction-tensor-criterion}
 \sum_{i=1}^d\Psi_i(\ell_i)
 =a\otimes t^{\mathrm{op}}-s\otimes b^{\mathrm{op}}.
\end{equation}
\item\label{item:delta-bounded-fraction-tensors}
There exist a fraction pair $Q=(s,a,b,t)$ for $u$ and tensors $\ell_1,\ldots,\ell_d\in\M\otimes_{\C}\M^{\mathrm{op}}$ satisfying \eqref{eq:all-fraction-tensor-criterion}.
\item\label{item:delta-finite-rank} There exist a fraction pair $Q=(s,a,b,t)$ for $u$ and a tuple $T=(T_1,\ldots,T_d)\in\F(H)^d$ such that
\begin{equation}
 \sum_{i=1}^d[\rho(X_i),T_i]=\Gamma[Q,P_\Omega].
 \label{eq:certificate}
\end{equation}
\end{enumerate}
In~\textup{(\ref{item:delta-tensors})}, $k_1,\ldots,k_d$ are uniquely determined by $u$.
The tensors $\ell_1,\ldots,\ell_d$ in \textup{(\ref{item:delta-all-fraction-tensors})} and \textup{(\ref{item:delta-bounded-fraction-tensors})} and the tuple $T$ in~\textup{(\ref{item:delta-finite-rank})} are uniquely determined by the fraction pair $Q$.
\end{theorem}

Here $P_\Omega$ denotes the orthogonal projection of $H$ onto $\C\Omega$.
By Definition~\ref{def:fraction-commutator},
\[
\Gamma[Q,P_\Omega]=sP_\Omega b-aP_\Omega t\in\F(H).
\]
The left-hand sides of equations~\eqref{eq:tensor-criterion} and~\eqref{eq:certificate} are values of degree-one Hochschild differentials with the corresponding coefficient bimodules; see Remark~\ref{rem:boundary-interpretation}.

\begin{remark}
When the equivalent assertions of Theorem~\ref{thm:delta} hold, the tensors in assertion~\textup{(\ref{item:delta-tensors})} are the free difference quotients $k_i=\partial_i u$ of $u$ with respect to $X_1,\ldots,X_d$;
see the implication $(1)\implies(2)$ in its proof.
The implication $(2)\implies(3)$ follows as $\Psi_i$ is a right multiplication map.
The implication $(3)\implies(4)$ is obtained by clearing denominators in $\ell_i$.
The implication $(4)\implies(5)$ follows by identifying tensors with finite-rank operators as in equations~\eqref{eq:theta-bounded} and~\eqref{eq:tensor-commutator-general}.
The key input for $(5)\implies(1)$ is that $\Div_{L^0(\M)}\A_X$ is a division ring, as guaranteed under $\Delta(X)=d$ by \cite[Theorem~1.1]{MSY23}.
\end{remark}

\subsection{Analytic criteria}
Condition~\textup{(\ref{item:delta-finite-rank})} in Theorem~\ref{thm:delta} requires a fraction pair $Q$ for which the commutator equation \eqref{eq:certificate} has the (unique) finite-rank solution $T\in\F(H)^d$.
The next theorem gives a rationality criterion that is valid for every fraction pair.
That is, \eqref{eq:certificate} admits a finite-rank solution either for every fraction pair or for none.

To formulate this analytic criterion, we use the dense subspace $\M\Omega:=\{x\Omega:x\in\M\}$ of $H$.
We use the following variant of the (liberation) dual systems in \cite[Definition~13.1]{Voi99}.
\begin{definition}\label{def:liberation-dual-system}
Let $X=(X_1,\ldots,X_d)\in\M^d$.
A \emph{liberation dual system} for $X$ is a tuple $(L_1,\ldots,L_d)$ of linear operators on $H$ such that
\[
 \M\Omega\subseteq\dom L_i\cap\dom L_i^*
 \qquad(1\leq i\leq d)
\]
and, on $\M\Omega$,
\[
 [X_j,L_i]=\delta_{ij}[X_i,P_\Omega]
 \qquad(1\leq i,j\leq d).
\]
We call this system \emph{normalized} if, on $\M\Omega$,
\begin{equation}\label{eq:liberation-normalization}
 \sum_{i=1}^dL_i=P_\Omega.
\end{equation}
\end{definition}
Proposition~\ref{prop:liberation-normalization} shows that every liberation dual system can be normalized.

\Needspace{6\baselineskip}
\begin{theorem}\label{thm:intrinsic-analytic}
Let $X=(X_1,\ldots,X_d)\in\M^d$, and suppose that $X_i-\lambda$ is invertible in $L^0(\M)$ for every $1\leq i\leq d$ and every $\lambda\in\C$.
Let $(L_1,\ldots,L_d)$ be a normalized liberation dual system for $X$ in the sense of Definition~\ref{def:liberation-dual-system}.
Then, for every $u\in L^0(\M)$, the following are equivalent:
\begin{enumerate}
\item $u\in\Div_{L^0(\M)}\A_X$.
\item\label{item:intrinsic-analytic-existence-criterion}  There exist a fraction pair $Q$ for $u$ and $T_1,\ldots,T_d\in\F(H)$ such that
\begin{equation}
 \sum_{i=1}^d[\rho(X_i),T_i]=\Gamma[Q,P_\Omega].
 \label{eq:intrinsic-analytic-criterion}
\end{equation}
\item\label{item:intrinsic-analytic-criterion} For every fraction pair $Q$ for $u$, there exist $T_1,\ldots,T_d\in\F(H)$ satisfying \eqref{eq:intrinsic-analytic-criterion}.
\item\label{item:intrinsic-analytic-fraction-commutators} For every fraction pair $Q$ for $u$, there exist $T_1,\ldots,T_d\in\F(H)$ such that
\begin{equation}\label{eq:intrinsic-analytic-coordinates}
 [\rho(X_i),T_i]=\overline{\Gamma[Q,L_i]}
 \qquad(1\leq i\leq d).
\end{equation}
\end{enumerate}
For each fixed fraction pair $Q$, the solution tuple $(T_1,\ldots,T_d)$ in assertions~\textup{(\ref{item:intrinsic-analytic-existence-criterion})}--\textup{(\ref{item:intrinsic-analytic-fraction-commutators})} is unique whenever it exists.
\end{theorem}

The new characterization in assertion~\textup{(\ref{item:intrinsic-analytic-fraction-commutators})} can be understood as follows.
Let $(T_1,\ldots,T_d)$ be a finite-rank solution of equation~\eqref{eq:intrinsic-analytic-coordinates} for a fixed fraction pair $Q$.
Summing over $i$ and using the normalization in equation~\eqref{eq:liberation-normalization}, we obtain
\[
 \sum_i[\rho(X_i),T_i]
 =\sum_i\overline{\Gamma[Q,L_i]}
 =\Gamma[Q,P_\Omega],
\]
so the same tuple solves equation~\eqref{eq:intrinsic-analytic-criterion}.
Conversely, we will see in Section~\ref{sec:liberation-dual-system-criterion}, after clearing denominators, the liberation relations ensure that the map $c\otimes e^{\mathrm{op}}\mapsto cL_i e$ isolates the $i$th summand, yielding equation~\eqref{eq:intrinsic-analytic-coordinates} for the given $T_i$ after cancellation and passage to the closure.

By Proposition~\ref{prop:liberation-delta}, the hypotheses on $X$ in Theorem~\ref{thm:intrinsic-analytic} imply $\Delta(X)=d$, the hypothesis of Theorem~\ref{thm:delta}.
These hypotheses are strictly stronger, as Example~\ref{ex:maximal-delta-without-liberation} in Section~\ref{sec:examples} exhibits a selfadjoint tuple with maximal $\Delta$ and regular scalar shifts that admits no liberation dual system.

Theorems~\ref{thm:delta} and~\ref{thm:intrinsic-analytic} allow the $X_i$ to be nonnormal.
When $X_1,\ldots,X_d$ are normal operators that generate $\M$ and admit a dual system $(D_1,\ldots,D_d)$, Theorem~\ref{thm:analytic} characterizes rationality by the finite rank of the fraction commutators $\Gamma[Q,D_i]$ for a given fraction pair $Q$.
If $u$ is rational in $X$, the proof identifies the common solution of equations~\eqref{eq:intrinsic-analytic-criterion} and~\eqref{eq:intrinsic-analytic-coordinates} for each fraction pair $Q$ as $T_i=\overline{\Gamma[Q,D_i]}$ for every $i$.

To state our third main result in a general setting, we use the following unbounded variant of the dual systems in \cite[Definition~5.1]{Voi98}.

\begin{definition}\label{def:dual-system}
Let $X=(X_1,\ldots,X_d)\in\M^d$.
A \emph{dual system} for $X$ is a tuple $(D_1,\ldots,D_d)$ of linear operators on $H$ such that
\begin{equation}\label{eq:analytic-hyp}
 \M\Omega\subseteq\dom D_i\cap\dom D_i^*,\qquad
 [X_j,D_i]|_{\M\Omega}=\delta_{ij}P_\Omega|_{\M\Omega}
 \qquad(1\leq i,j\leq d).
\end{equation}
\end{definition}

For a tuple $X=(X_1,\ldots,X_d)\in\M^d$, write $W^*(X)=W^*(X_1,\ldots,X_d)$ for the unital von Neumann subalgebra of $\M$ generated by $X_1,\dots,X_d$.

\Needspace{6\baselineskip}
\begin{theorem}\label{thm:analytic}
Let $X$ be a tuple of bounded normal operators and let $(D_1,\ldots,D_d)$ be a dual system for $X$ on $L^2(W^*(X),\tau)$.
For every $u\in L^0(W^*(X))$ and every fraction pair $Q=(s,a,b,t)\in W^*(X)^4$ for $u$, the following conditions are equivalent.
\begin{enumerate}
\item $u\in\Div_{L^0(W^*(X))}\A_X$.
\item For every $i$, $\Gamma[Q,D_i]$ has finite rank.
\item\label{item:analytic-fraction-commutators} For every $i$, $\overline{\Gamma[Q,D_i]}\in\F\bigl(L^2(W^*(X),\tau)\bigr)$.
\end{enumerate}
When these equivalent conditions hold, the rank of each $\overline{\Gamma[Q,D_i]}$ is independent of $Q$.
\end{theorem}

Proposition~\ref{prop:fraction-commutator-rank} computes these ranks for $u=r(X)$, where $r$ is a nonzero rational function.
Let $A=A_0+\sum_{j=1}^d A_jx_j$ be the matrix in a minimal linear representation of $r$.
For every fraction pair $Q$ for $u$, Proposition~\ref{prop:fraction-commutator-rank} gives
\[
 \rank\overline{\Gamma[Q,D_i]}
 =\rank_{\C} A_i
 =\text{the tensor rank of }\partial_i r
 \qquad(1\leq i\leq d),
\]
where $\partial_i r$ is the $i$th free difference quotient of $r$.
In one variable, this tensor rank equals the Hankel rank of the associated Laurent series, whose rationality is characterized by finiteness of that rank in Kronecker's theorem \cite{Kro81}.
Theorem~\ref{thm:analytic} thus gives a noncommutative counterpart to Kronecker's criterion.

Even among normal generators, the hypotheses of Theorem~\ref{thm:analytic} are strictly stronger than those of Theorem~\ref{thm:intrinsic-analytic}.
The arcsine generator $X$ in Example~\ref{ex:liberation-without-duals} admits a normalized liberation dual system, and $X-\lambda$ is regular for every $\lambda\in\C$.
The same generator admits no dual system.

Theorem~\ref{thm:analytic} gives the rationality criterion for free Haar unitaries \cite{DR97,Lin00} and for free semicircular operators \cite{Miy22}.
In both cases the generators are normal and admit bounded dual systems, for which $\Gamma[Q,D_i]$ is bounded on $H$ and equals its closure.

For $\M=W^*(X)$, the equality
\[
 \Div_{\M}\A_X=\Div_{L^0(\M)}\A_X\cap\M
\]
holds when $X$ is a tuple of free Haar unitaries \cite[Theorem~1.1]{Lin00} or a standard free semicircular tuple \cite[Theorem~4.5]{Miy22}.
Subsection~\ref{subsec:bounded-division-closures} shows that this equality can fail even for a selfadjoint tuple admitting a bounded dual system.

\subsection{Organization}\label{subsec:organization}

Section~\ref{sec:preliminaries} collects the algebraic preliminaries on division closures and the free field, followed by fractions of affiliated operators, finite-rank tensor realizations and the consequences of maximal $\Delta$.
Section~\ref{sec:maximal-delta} proves Theorem~\ref{thm:delta} using denominator clearing and the tensor realization of finite-rank operators.
Section~\ref{sec:liberation-dual-system-criterion} proves that liberation dual systems give maximal $\Delta$ when every $X_i-\lambda$ is regular, and establishes Theorem~\ref{thm:intrinsic-analytic} by recovering algebraic $L^2$ tensors.
Section~\ref{sec:dual-system-criterion} constructs normalized liberation dual systems from dual systems.
It then deduces Theorem~\ref{thm:analytic} directly from Theorem~\ref{thm:intrinsic-analytic}.
It also computes the ranks of the fraction commutator closures from minimal linear representations and identifies the rank with Hankel rank when $d=1$.
Section~\ref{sec:examples} separates maximal $\Delta$ from the existence of a liberation dual system and distinguishes the two analytic hypotheses.
A circular example illustrates the role of normality in Theorem~\ref{thm:analytic}.
The final example shows that the division closure in $\M$ can be strictly smaller than the bounded part of the division closure in $L^0(\M)$, even for selfadjoint generators with a bounded dual system.

\subsection*{Acknowledgments}
The author is grateful to R.~Speicher for sharing with him the question of rationality for affiliated operators during the COVID-19 quarantine, and to R.~Speicher and T.~Mai for stimulating discussions during that period.
The author also thanks A.~Miyagawa for his feedback after reading a preliminary version of this manuscript.

\subsection*{Disclosure of AI Assistance}
The author used \href{https://github.com/vtejdn/acme}{ACME}, a framework for autonomous and collaborative mathematical exploration, in this research.
The author set the research objectives, contributed mathematical insights, and guided key decisions.
Within this framework, AI agents autonomously explored possible solutions and proof strategies.

The research materials and exploration records are publicly available in the author's \href{https://github.com/vtejdn/yin-math-research}{GitHub repository} under Task F0001 (versions 1--7).
All arguments and conclusions presented in this manuscript have been carefully verified by the author.
The author takes responsibility for the content of this manuscript.

\section{Preliminaries}\label{sec:preliminaries}

We first recall division closures, the free field and an algebraic tensor criterion for membership in a division subring.
We then discuss fractions of affiliated operators, realize algebraic tensors as finite-rank operators and recall the consequences of maximal $\Delta$.

\subsection{Division closures, the free field and linearization}\label{subsec:free-field}

\begin{definition}\label{def:rationality}
For a unital subalgebra $A$ of a unital algebra $B$, its \emph{division closure} $\Div_B A$ is the smallest unital subalgebra containing $A$ and closed under inversion of its elements that are invertible in $B$; see \cite[Exercise~7.1.4]{Coh85} and \cite[Definition~4.10]{MSY23}.
An affiliated operator $u$ is \emph{rational in $X$} if $u\in\Div_{L^0(\M)}\A_X$.
\end{definition}

The related notion of \emph{rational closure} is defined using entries of matrix inverses \cite[Definition~4.7]{MSY23}; see also \cite[Section~7]{Coh71}.
It contains the division closure and may be strictly larger in general \cite[Exercise~7.1.4]{Coh85}.
Since the division and rational closures coincide when $\Delta(X)=d$ (see Subsection~\ref{subsec:delta-tensor-maps}), we use \emph{rational in $X$} for membership in the division closure.

The \emph{free field} $\freefield{x_1,\ldots,x_d}$ is the universal division ring of fractions of the free polynomial algebra $\C\langle x_1,\ldots,x_d\rangle$; see \cite[Section~4.1 and Theorem~4.1]{MSY23} and \cite[Corollary~7.5.14]{Coh06}.
A square matrix $A$ of size $n$ over the polynomial algebra is \emph{full} if every factorization $A=BC$ over that algebra, with $B$ of size $n\times k$ and $C$ of size $k\times n$ for some $k\geq 0$, satisfies $k\geq n$ \cite[p.~613]{Coh71}.
The free field is the universal localization of the polynomial algebra at all full square matrices \cite[Theorem~7.5.13 and Corollary~7.5.14]{Coh06}.
We will call elements of the free field \emph{(noncommutative) rational functions}.

The free difference quotients on noncommutative polynomials \cite[Section~3, pp.~192--193]{Voi98} extend uniquely to the free field \cite[Section~2.1, (4), pp.~298--299]{BD78}.
An explicit construction of these extensions on noncommutative rational functions is given in \cite[Section~4.4, Definition~4.4 and Corollary~4.11]{KV12}.
For $1\leq i\leq d$, we use the linear identification $u\otimes v\mapsto u\otimes v^{\mathrm{op}}$ to write the extended maps as
\[
 \partial_i:\freefield{x_1,\ldots,x_d}\longrightarrow
 \freefield{x_1,\ldots,x_d}\otimes_{\C}
 \freefield{x_1,\ldots,x_d}^{\mathrm{op}}.
\]
For $c\in\C$ and $u,v\in\freefield{x_1,\ldots,x_d}$, with $u\ne0$, they satisfy
\begin{equation}\label{eq:free-field-difference-rules}
\begin{aligned}
 \partial_i c&=0,& \partial_i x_j&=\delta_{ij}1\otimes1,\\
 \partial_i(uv)&=(1\otimes v^{\mathrm{op}})\partial_i u
                  +(u\otimes1)\partial_i v,&
 \partial_i(u^{-1})&=-(u^{-1}\otimes(u^{-1})^{\mathrm{op}})\partial_i u.
\end{aligned}
\end{equation}

\begin{definition}\label{def:linear-representation}
Let $r\in\freefield{x_1,\ldots,x_d}$ be a rational function.
A \emph{linear representation} of $r$ is a tuple $(u,A,v)$ consisting of a full linear matrix
\[
 A=A_0+\sum_{j=1}^d A_jx_j,\qquad A_0,\ldots,A_d\in M_n(\C),\qquad n\geq 1,
\]
a row vector $u\in M_{1,n}(\C)$, and a column vector $v\in M_{n,1}(\C)$ such that $r=uA^{-1}v$ in the free field.
The \emph{dimension} of this representation is $n$.
The representation is \emph{minimal} if $n$ is the least possible dimension among all representations of $r$ \cite[Section~3.1]{CR94}.
\end{definition}

\begin{theorem}[{\cite[Section~5.8 and Theorem~7.1.2]{Coh85}}]\label{thm:linearization}
Each rational function $r\in\freefield{x_1,\ldots,x_d}$ admits a linear representation in the sense of Definition~\ref{def:linear-representation}, and hence a minimal one.
\end{theorem}

\begin{theorem}[{\cite[Section~4.2, Proposition~7]{CR94}}]\label{thm:minimal-linear-representation}
Let $(u,A,v)$ be a linear representation of a nonzero rational function.
It is minimal if and only if the entries of $uA^{-1}$ are linearly independent over $\C$ and so are the entries of $A^{-1}v$.
\end{theorem}

\subsection{Coordinate maps over division subrings}

For a division subring $F$ in a unital algebra $U$, the proof below constructs, for $u\notin F$, a right $F$-linear coordinate map $\phi_u:U\longrightarrow F$ satisfying $\phi_u(1)=0$ and $\phi_u(u)=1$.
This map distinguishes $u\otimes_F1$ from $1\otimes_Fu$, which is recorded as the lemma below.
This is the key lemma used in the proof of Theorem~\ref{thm:delta} to force an affiliated operator to be a rational function in $X$.

\begin{lemma}\label{lem:equalizer}
Let $U$ be a nonzero unital algebra and let $F\subseteq U$ be a division subring containing the identity of $U$.
Then, for every $u\in U$,
\[
 u\otimes_F1=1\otimes_Fu\text{ in }U\otimes_FU
 \quad\Longleftrightarrow\quad u\in F.
\]
\end{lemma}
\begin{proof}
If $u\in F$, the identity follows from the balancing relation
\[
 (v_1a)\otimes_F v_2=v_1\otimes_F (av_2),\qquad v_1,v_2\in U,\quad a\in F.
\]
If $u\notin F$, the elements $1,u$ are right $F$-linearly independent: $v_1+uv_2=0$ with $v_2\ne0$ would give $u=-v_1v_2^{-1}\in F$.
Extend $(1,u)$ to a basis of the right $F$-module $U$ by \cite[Theorem~IV.2.4]{Hun74}.
The coordinate map $\phi_u$ at $u$ is right $F$-linear (see \cite[Theorem~IV.2.1 \& Theorem~4.11]{Hun74}) and satisfies $\phi_u(1)=0$, $\phi_u(u)=1$.
The map $(v_1,v_2)\mapsto\phi_u(v_1)v_2$ is middle linear and hence induces an additive map \cite[Definition~IV.5.1 and Theorem~IV.5.2]{Hun74}
\[
 U\otimes_FU\longrightarrow U,\qquad v_1\otimes_Fv_2\longmapsto\phi_u(v_1)v_2
\]
which sends $u\otimes_F1-1\otimes_Fu$ to $1$, so this tensor cannot vanish.
\end{proof}

\subsection{Fractions of affiliated operators}\label{subsec:rational-affiliated}

We now turn to fraction representations in $L^0(\M)$.
Affiliation means commutation with each unitary in the commutant, including invariance of the domain; see \cite[Section 7.2]{AP18} and \cite[Section 5.1]{MSY23}.

\begin{definition}\label{def:regular}
An element $s\in\M$ is \emph{regular} if it is invertible in $L^0(\M)$.
\end{definition}

Regularity is equivalent to injectivity and dense range on $H$ \cite[Lemmas 5.1--5.2]{MSY23}.
For $u\in L^0(\M)$, a fraction pair is given by
\begin{equation}
 s:=(1+uu^*)^{-1/2},\quad a:=su,\qquad
 t:=(1+u^*u)^{-1/2},\quad b:=ut,
 \label{eq:canonical}
\end{equation}
where $s,a,b,t\in\M$, $s,t$ are regular, and $u=s^{-1}a=bt^{-1}$.
These facts follow from functional calculus and polar decomposition for affiliated operators.

\begin{definition}\label{def:fraction-terminology}
For a fraction pair $Q=(s,a,b,t)$ of $u$ as in Definition~\ref{def:fraction-pair}, the representations $s^{-1}a$ and $bt^{-1}$ are called a \emph{left fraction} and a \emph{right fraction}, respectively; see \cite[Section~8.2.1]{Luc02}.

For $u_1,\ldots,u_m\in L^0(\M)$, a regular $s\in\M$ is a \emph{common left denominator} if $su_j\in\M$ for every $j$; a regular $t\in\M$ is a \emph{common right denominator} if $u_jt\in\M$ for every $j$; see \cite[Section~1]{Lin00} and \cite[Lemma~4.1]{Miy22}.
\end{definition}

Common denominators for finite families are used in \cite[Section~1, p.~2 of the arXiv version]{Lin00}.
The following proof iterates the two-element construction of \cite[Lemma~4.1]{Miy22}.

\begin{lemma}\label{lem:clearing}
For finite families $u_1,\ldots,u_m$ and $v_1,\ldots,v_n$ in $L^0(\M)$ there are regular $s,t\in\M$ such that
\[
 su_j\in\M\quad(1\leq j\leq m),\qquad
 v_kt\in\M\quad(1\leq k\leq n).
\]
\end{lemma}
\begin{proof}
Choose a regular $s\in\M$ with $su_1\in\M$ using \eqref{eq:canonical}.
Suppose $su_j\in\M$ for $1\leq j\leq\ell$, where $1\leq\ell<m$.
Using \eqref{eq:canonical}, write
\[
 u_{\ell+1}=s_1^{-1}a_1,\qquad ss_1^{-1}=s_2^{-1}a_2,
\]
where $s_1,s_2,a_1,a_2\in\M$ and $s_1,s_2$ are regular.
Then $s':=s_2s=a_2s_1$ is regular, and
\[
 s'u_j=s_2(su_j)\in\M\quad(1\leq j\leq\ell),\qquad
 s'u_{\ell+1}=a_2a_1\in\M.
\]
This completes the induction for the left denominators.
The proof for right denominators is similar.
\end{proof}

\begin{lemma}[{\cite[proof of Lemma~4.1]{Miy22}}]\label{lem:common}
If $s_1,s_2\in\M$ are regular, they have a common left multiple $xs_1=ys_2$ with $x,y\in\M$ regular.
If $t_1,t_2\in\M$ are regular, they have a common right multiple $t_1x'=t_2y'$ with $x',y'\in\M$ regular.
Both common multiples are regular.
\end{lemma}
\begin{proof}
Choose a left fraction $s_1s_2^{-1}=(s')^{-1}a$ using \eqref{eq:canonical}.
Thus we may take $x:=s'$ and $y:=a$.
Indeed, $s's_1=as_2$, and $a=s's_1s_2^{-1}$ is invertible in $L^0(\M)$, hence regular.
A right fraction $t_1^{-1}t_2=b_3t_3^{-1}$ gives the second identity, with $b_3=t_1^{-1}t_2t_3$ regular.
Products of regular elements are regular.
\end{proof}

Linnell \cite[Theorem~1.1]{Lin00} characterizes rational affiliated operators over a free group by finite-rank conditions on bounded operators obtained from left and right fractions.
Miyagawa \cite[Section~4]{Miy22} gives a semicircular analogue using the expressions
\[
 s r_i^*b-a r_i^*t=\Gamma[Q,r_i^*],
\]
where $Q=(s,a,b,t)$ is a fraction pair and $r_i^*$ are the bounded right annihilation operators.

\subsection{Algebraic tensors and finite-rank operators}\label{subsec:finite-rank-tensors}

Recall that $\F(H)$ is the space of bounded finite-rank operators on $H$ and that $\rho(x)$ denotes right multiplication by $x\in\M$ on $H$.
We use the commutator convention $[S,T]:=ST-TS$.

Let $S_2(H)$ be the Schatten class of Hilbert--Schmidt operators, that is, the Hilbert space of operators $S\in\B(H)$ with $\operatorname{Tr}(S^*S)<\infty$, equipped with the Hilbert--Schmidt norm $\|S\|_{\operatorname{HS}}:=\operatorname{Tr}(S^*S)^{1/2}$.
Here $\operatorname{Tr}$ is the usual trace on $\B(H)$ \cite[Section 8.4]{AP18}.

Let $J$ be Tomita's conjugation operator on $H$, given by $J(x\Omega):=x^*\Omega$; see \cite[Section~7.1.1]{AP18} and \cite[p.~13]{MSY23}.
It satisfies $J^2=I$ and $\rho(x)=Jx^*J$.
Following \cite[Section~1.1, equation~(1.1)]{CS05}, define
\begin{equation}
 \Theta(u\otimes v)\xi:=u\langle Jv,\xi\rangle,
 \qquad \Theta:H\otimes_{\C}H\longrightarrow\F(H).
 \label{eq:theta}
\end{equation}
Here $u,v\in H=L^2(\M,\tau)\subset L^0(\M)$ and $\xi\in H$.
The map $\Theta$ is a linear bijection and extends to a unitary
\[
 \Theta:H\bar\otimes H\longrightarrow S_2(H).
\]
The actions in \cite[equation~(1.2)]{CS05} give, for $a,b\in\M$ and $u,v\in H$,
\[
 \Theta(au\otimes vb)=a\Theta(u\otimes v)b.
\]

Via the linear identification $u\otimes v\mapsto u\otimes v^{\mathrm{op}}$, we also regard $\Theta$ as a map on the corresponding subspace of $L^0(\M)\otimes_{\C}L^0(\M)^{\mathrm{op}}$.
For $u,v\in\M$, equation~\eqref{eq:theta} reads
\begin{equation}\label{eq:theta-bounded}
 \Theta(u\otimes v^{\mathrm{op}})=uP_\Omega v.
\end{equation}
For $x\in\M$ and $u,v\in H$, the right multiplication operator $\rho(x)$ satisfies
\begin{align*}
 \rho(x)\Theta(u\otimes v^{\mathrm{op}})
   &=\Theta(ux\otimes v^{\mathrm{op}}),\\
 \Theta(u\otimes v^{\mathrm{op}})\rho(x)
   &=\Theta(u\otimes(xv)^{\mathrm{op}}).
\end{align*}
Subtracting gives
\begin{equation}
 [\rho(x),\Theta(u\otimes v^{\mathrm{op}})]
   =\Theta(ux\otimes v^{\mathrm{op}}-u\otimes(xv)^{\mathrm{op}}).
 \label{eq:tensor-commutator}
\end{equation}
Extending equation~\eqref{eq:tensor-commutator} by linearity to finite sums gives
\begin{equation}
 [\rho(x),\Theta(k)]=\Theta(\Psi_x(k)),
 \qquad k\in H\otimes_{\C}H.
 \label{eq:tensor-commutator-general}
\end{equation}
Here $\Psi_x$ is the right multiplication map of Definition~\ref{def:psi}.

\subsection{The invariant \texorpdfstring{$\Delta$}{Delta} and free-field realization}\label{subsec:delta-tensor-maps}

\Needspace{20\baselineskip}
\begin{definition}\label{def:delta}
The von Neumann tensor product $\M\bar\otimes\M^{\mathrm{op}}$, where $\M^{\mathrm{op}}$ is the opposite algebra, acts on $S_2(H)$ by
\[
 (a\otimes b^{\mathrm{op}})\cdot S:=aSb,\qquad a,b\in\M.
\]
For $X=(X_1,\ldots,X_d)\in\M^d$, set
\[
 \mathcal C_X:=
 \left\{T=(T_1,\ldots,T_d)\in\F(H)^d:
 \sum_{i=1}^d[\rho(X_i),T_i]=0\right\}.
\]
We call an element $T\in\mathcal C_X$ a \emph{finite-rank cycle}.
Such tuples represent Hochschild $1$-cycles with Hilbert--Schmidt coefficients; see \cite{CS05}.
Following \cite{CS05}, define
\begin{equation}
 \Delta(X):=d-\dim_{\M\bar\otimes\M^{\mathrm{op}}}
 \overline{\mathcal C_X}^{\operatorname{HS}}.
 \label{eq:delta-definition}
\end{equation}
The closure is taken in $S_2(H)^d$.
The dimension is the von Neumann dimension \cite[Definition~1.10]{Luc02} for the product trace $(a\otimes b^{\mathrm{op}})\mapsto\tau(a)\tau(b)$; see also \cite[Definition~2.10]{CS05}.
Thus $\Delta(X)\leq d$.
\end{definition}

We recall the finite-rank characterization of maximal $\Delta$ and the resulting free-field realization of the division closure, both from \cite{MSY23}.

\begin{remark}\label{rem:cycle-characterization}
By \cite[Section~3 and Remark~3.2]{MSY23}, for every $X\in\M^d$,
\[
 \Delta(X)=d
 \quad\Longleftrightarrow\quad \mathcal C_X=\{0\}
 \quad\Longleftrightarrow\quad
 \left\{T\in\F(H)^d:\sum_i[T_i,X_i]=0\right\}=\{0\}.
\]
\end{remark}

\begin{lemma}[{\cite[Theorems 1.1, 3.1 and 5.6]{MSY23}}]
\label{lem:free-field}
If $\Delta(X)=d$, evaluation at $X$ extends to an embedding of the free field $\freefield{x_1,\ldots,x_d}$ into $L^0(\M)$.
Its image is the division subring $\Div_{L^0(\M)}\A_X$.
\end{lemma}

If $\Delta(X)=d$, the rational closure and division closure of $\A_X$ in $L^0(\M)$ coincide \cite[Theorem~1.1 and Remark~5.10]{MSY23}.

We end this subsection with a remark that the above lemma can be extended to noncommutative meromorphic functions; see a recent preprint \cite{AMTS26} by Augat, Martin, Tornes and Shamovich.
They extend evaluation at $X$ from the ring of noncommutative power series to an embedding into $L^0(\M)$ under a similar maximality condition on the commutator equation $\sum_i[T_i,X_i]=0$.
Instead of finite rank operators, they consider the solutions of the commutator equation among tuples of trace-class operators.
They also give a selfadjoint triple with maximal $\Delta$ for which this commutator equation has a nonzero trace-class solution \cite[Remark~7.5]{AMTS26}.
This distinction recalls the classical commutator with the Hardy projection, which has finite rank precisely for rational functions and is of trace class precisely for the Besov space $B^1_{1,1}(\mathbb T)$ by Peller's theorem \cite[Theorem~5.3]{Pel24}.
\section{Rationality under maximal \texorpdfstring{$\Delta$}{Delta}}\label{sec:maximal-delta}

We prove Theorem~\ref{thm:delta} in this section.
By Definition~\ref{def:psi}, $\Psi_i$ is right multiplication by $X_i\otimes1-1\otimes X_i^{\mathrm{op}}$.
In particular,
\begin{equation}\label{eq:psi-left-linearity}
 \Psi_i(rk)=r\Psi_i(k),\qquad
 r,k\in L^0(\M)\otimes_{\C}L^0(\M)^{\mathrm{op}}.
\end{equation}

\begin{lemma}\label{lem:tensor-clearing}
Given $k_1,\ldots,k_d\in L^0(\M)\otimes_{\C}L^0(\M)^{\mathrm{op}}$, there are regular $s,t\in\M$ such that
\[
 (s\otimes t^{\mathrm{op}})k_i
 \in\M\otimes_{\C}\M^{\mathrm{op}}\qquad(1\leq i\leq d).
\]
\end{lemma}
\begin{proof}
For simplicity of notation, we give the calculation for simple tensors
\[
 k_i=u_i\otimes v_i^{\mathrm{op}},\qquad u_i,v_i\in L^0(\M).
\]
The same argument applies to general finite sums.
Applying Lemma~\ref{lem:clearing} to $u_1,\ldots,u_d$ gives a regular $s\in\M$ such that
\[
 su_i\in\M\qquad(1\leq i\leq d).
\]
Applying its right-denominator assertion to $v_1,\ldots,v_d$ gives a regular $t\in\M$ such that
\[
 v_it\in\M\qquad(1\leq i\leq d).
\]
Then
\[
 (s\otimes t^{\mathrm{op}})k_i
 =su_i\otimes(v_it)^{\mathrm{op}}
 \in\M\otimes_{\C}\M^{\mathrm{op}}.
\]
\end{proof}

\begin{lemma}\label{lem:tensor-injectivity}
Suppose $\Delta(X)=d$.
The map
\[
 (L^0(\M)\otimes_{\C}L^0(\M)^{\mathrm{op}})^d\longrightarrow
 L^0(\M)\otimes_{\C}L^0(\M)^{\mathrm{op}},\qquad
 (k_i)\longmapsto\sum_i\Psi_i(k_i)
\]
is injective.
\end{lemma}
\begin{proof}
Suppose $\sum_i\Psi_i(k_i)=0$.
Choose regular $s,t\in\M$ as in Lemma~\ref{lem:tensor-clearing} and set $T_i:=\Theta((s\otimes t^{\mathrm{op}})k_i)\in\F(H)$.
By \eqref{eq:tensor-commutator-general} and \eqref{eq:psi-left-linearity},
\[
 \sum_i[\rho(X_i),T_i]
 =\Theta\left((s\otimes t^{\mathrm{op}})\sum_i\Psi_i(k_i)\right)=0.
\]
Remark~\ref{rem:cycle-characterization} gives $T_i=0$ for every $i$.
Injectivity of $\Theta$ implies $(s\otimes t^{\mathrm{op}})k_i=0$.
Left multiplication by $s^{-1}\otimes(t^{-1})^{\mathrm{op}}$ then gives $k_i=0$.
\end{proof}

\begin{lemma}\label{lem:tensor-certificate}
Let $Q=(s,a,b,t)$ be a fraction pair for $u\in L^0(\M)$, and let $\ell_1,\ldots,\ell_d\in H\otimes_{\C}H$.
Set $T_i:=-\Theta(\ell_i)\in\F(H)$.
Then
\[
 \sum_i\Psi_i(\ell_i)
 =a\otimes t^{\mathrm{op}}-s\otimes b^{\mathrm{op}}
\]
if and only if
\[
 \sum_i[\rho(X_i),T_i]=\Gamma[Q,P_\Omega].
\]
\end{lemma}
\begin{proof}
By \eqref{eq:tensor-commutator-general} and \eqref{eq:theta-bounded},
\[
 \sum_i[\rho(X_i),T_i]
 =-\Theta\left(\sum_i\Psi_i(\ell_i)\right),
 \qquad
 \Gamma[Q,P_\Omega]
 =-\Theta(a\otimes t^{\mathrm{op}}-s\otimes b^{\mathrm{op}}).
\]
The equivalence follows from the injectivity of $\Theta$.
\end{proof}

By Lemma~\ref{lem:free-field}, $\Div_{L^0(\M)}\A_X$ is a division subring of $L^0(\M)$ whenever $\Delta(X)=d$.
We prove Theorem~\ref{thm:delta} through $(1)\implies(2)\implies(3)\implies(4)\implies(5)\implies(1)$.

\begin{proof}[Proof of Theorem~\ref{thm:delta}]
In this proof, set $U:=L^0(\M)$ and $F:=\Div_U\A_X$.

$(1)\implies(2)$.
By \cite[Theorem~1.1, (i) and (v)]{MSY23}, the assumption $\Delta(X)=d$ implies that $X_1,\ldots,X_d$ generate the free field: evaluation identifies $\freefield{x_1,\ldots,x_d}$ with $F$.
The free difference quotients from Subsection~\ref{subsec:free-field}, transported through this evaluation isomorphism, define linear maps
\[
 \partial_i:F\longrightarrow F\otimes_{\C}F^{\mathrm{op}}.
\]
They satisfy the rules in equation~\eqref{eq:free-field-difference-rules} with $x_j$ replaced by $X_j$.
Put $\delta v:=v\otimes1-1\otimes v^{\mathrm{op}}$.
Direct expansion gives
\[
 \delta(v_1v_2)=(1\otimes v_2^{\mathrm{op}})\delta v_1+(v_1\otimes1)\delta v_2,
 \qquad
 \delta(v^{-1})=-(v^{-1}\otimes(v^{-1})^{\mathrm{op}})\delta v.
\]
For constants, $\delta c=0$, while for generators, $\delta X_j=\Psi_j(1\otimes1)=\sum_i\Psi_i(\partial_iX_j)$.
We now show that this identity extends from constants and generators to all rational elements using the product and inverse rules.

If $\delta v_\ell=\sum_i\Psi_i(\partial_i v_\ell)$ for $\ell=1,2$, the product rules and equation~\eqref{eq:psi-left-linearity} give
\begin{align*}
 \delta(v_1v_2)
 &=(1\otimes v_2^{\mathrm{op}})\sum_i\Psi_i(\partial_i v_1)
   +(v_1\otimes1)\sum_i\Psi_i(\partial_i v_2)\\
 &=\sum_i\Psi_i\bigl((1\otimes v_2^{\mathrm{op}})\partial_i v_1
                  +(v_1\otimes1)\partial_i v_2\bigr)\\
 &=\sum_i\Psi_i\bigl(\partial_i(v_1v_2)\bigr).
\end{align*}
Similarly, if $0\ne v\in F$ satisfies $\delta v=\sum_i\Psi_i(\partial_i v)$, the inverse rules give
\begin{align*}
 \delta(v^{-1})
 &=-(v^{-1}\otimes(v^{-1})^{\mathrm{op}})\sum_i\Psi_i(\partial_i v)\\
 &=\sum_i\Psi_i\bigl(-(v^{-1}\otimes(v^{-1})^{\mathrm{op}})\partial_i v\bigr)\\
 &=\sum_i\Psi_i\bigl(\partial_i(v^{-1})\bigr).
\end{align*}
Together with linearity, this shows that the elements satisfying the identity form a division subalgebra of $F$ containing $\C$ and $X_1,\ldots,X_d$.
Since these generate $F$, the identity holds for every rational element $u\in F$.
This proves equation~\eqref{eq:tensor-criterion} with $k_i:=\partial_i u$; thus assertion~\textup{(\ref{item:delta-tensors})} holds.
Lemma~\ref{lem:tensor-injectivity} proves uniqueness of the tensors in assertion~\textup{(\ref{item:delta-tensors})}.

$(2)\implies(3)$.
Given tensors $k_i$ as in assertion~\textup{(2)}, for any fraction pair $Q=(s,a,b,t)$ set $\ell_i:=(s\otimes t^{\mathrm{op}})k_i$.
By equation~\eqref{eq:psi-left-linearity},
\[
 \sum_i\Psi_i(\ell_i)
 =(s\otimes t^{\mathrm{op}})
   (u\otimes1-1\otimes u^{\mathrm{op}})
 =a\otimes t^{\mathrm{op}}-s\otimes b^{\mathrm{op}},
\]
so assertion~\textup{(\ref{item:delta-all-fraction-tensors})} holds.

$(3)\implies(4)$.
Choose a fraction pair $Q=(s,a,b,t)$ for $u$ and corresponding tensors $\ell_i$ in assertion~\textup{(3)}.
By Lemma~\ref{lem:tensor-clearing}, there are regular $s',t'\in\M$ such that
\[
 (s'\otimes(t')^{\mathrm{op}})\ell_i
 \in\M\otimes_{\C}\M^{\mathrm{op}}.
\]
By equation~\eqref{eq:psi-left-linearity} and equation~\eqref{eq:all-fraction-tensor-criterion},
\[
 \sum_i\Psi_i\bigl((s'\otimes(t')^{\mathrm{op}})\ell_i\bigr)
 =s'a\otimes(tt')^{\mathrm{op}}-s's\otimes(bt')^{\mathrm{op}}.
\]
Since $(s's,s'a,bt',tt')$ is a fraction pair for $u$, assertion~\textup{(\ref{item:delta-bounded-fraction-tensors})} holds.

For fixed $Q$, uniqueness in assertions \textup{(\ref{item:delta-all-fraction-tensors})} and \textup{(\ref{item:delta-bounded-fraction-tensors})} follows from Lemma~\ref{lem:tensor-injectivity}, applied to the difference of two solutions of equation~\eqref{eq:all-fraction-tensor-criterion}.

$(4)\implies(5)$.
Let $Q$ and $\ell_i$ be as in assertion~\textup{(\ref{item:delta-bounded-fraction-tensors})}.
Lemma~\ref{lem:tensor-certificate} gives the finite-rank operators $T_i:=-\Theta(\ell_i)$ satisfying equation~\eqref{eq:certificate}.
Thus assertion~\textup{(\ref{item:delta-finite-rank})} holds.

$(5)\implies(1)$.
Suppose assertion~\textup{(\ref{item:delta-finite-rank})} holds.
Let $Q=(s,a,b,t)$ be the fraction pair and $T_1,\ldots,T_d\in\F(H)$ the operators satisfying equation~\eqref{eq:certificate}.
Via $\Theta$, define a linear map $\pi:\F(H)\longrightarrow U\otimes_FU$ by
\[
 \pi(\Theta(v_1\otimes v_2^{\mathrm{op}})):=v_1\otimes_Fv_2,
 \qquad v_1,v_2\in H.
\]
Note that $\Theta$ is a linear bijection, so $\pi$ is well defined.
For a rank-one operator $\Theta(v_1\otimes v_2^{\mathrm{op}})$, equation~\eqref{eq:tensor-commutator-general} gives
\begin{align*}
 \pi\bigl([\rho(X_i),\Theta(v_1\otimes v_2^{\mathrm{op}})]\bigr)
 &=\pi\bigl(\Theta(\Psi_i(v_1\otimes v_2^{\mathrm{op}}))\bigr)\\
 &=(v_1X_i)\otimes_Fv_2-v_1\otimes_F(X_iv_2)\\
 &=0.
\end{align*}
The last equality uses $X_i\in F$ and the balancing relation
\[
 (uf)\otimes_Fw=u\otimes_F(fw),\qquad u,w\in U,\quad f\in F.
\]
Every finite-rank operator is a finite sum of rank-one operators, so linearity gives $\pi([\rho(X_i),T])=0$ for every $T\in\F(H)$.

Applying $\pi$ to equation~\eqref{eq:certificate} and using $\Theta(v_1\otimes v_2^{\mathrm{op}})=v_1P_\Omega v_2$ for $v_1,v_2\in\M$ (see equation~\eqref{eq:theta-bounded}), we obtain
\begin{align*}
 0&=\pi\left(\sum_i[\rho(X_i),T_i]\right)\\
  &=\pi(\Gamma[Q,P_\Omega])\\
  &=\pi(sP_\Omega b-aP_\Omega t)\\
  &=\pi\bigl(\Theta(s\otimes b^{\mathrm{op}})-\Theta(a\otimes t^{\mathrm{op}})\bigr)\\
  &=s\otimes_Fb-a\otimes_Ft.
\end{align*}
Hence $a\otimes_Ft=s\otimes_Fb$.

To clear the denominators, consider the map
\[
 \Phi_{s,t}:U\otimes_FU\longrightarrow U\otimes_FU,\qquad
 \Phi_{s,t}(u\otimes_Fw):=s^{-1}u\otimes_Fwt^{-1}.
\]
It is well defined because it preserves the balancing relation.
Applying $\Phi_{s,t}$ to $a\otimes_Ft=s\otimes_Fb$ yields
\[
 s^{-1}a\otimes_Ftt^{-1}=s^{-1}s\otimes_Fbt^{-1}.
\]
Since $u=s^{-1}a=bt^{-1}$, this is
\[
 u\otimes_F1=1\otimes_Fu.
\]
Lemma~\ref{lem:equalizer} therefore gives $u\in F=\Div_U\A_X$, so assertion~\textup{(1)} holds.

For fixed $Q$, let $(T_1,\ldots,T_d)$ and $(S_1,\ldots,S_d)$ be two finite-rank solutions of equation~\eqref{eq:certificate}.
Their right-hand sides are both $\Gamma[Q,P_\Omega]$, so subtraction gives
\begin{align*}
 \sum_i[\rho(X_i),T_i-S_i]
 &=\sum_i[\rho(X_i),T_i]-\sum_i[\rho(X_i),S_i]=0.
\end{align*}
Each $T_i-S_i$ has finite rank.
Since $\Delta(X)=d$, Remark~\ref{rem:cycle-characterization} implies $T_i-S_i=0$ for every $i$.
Hence the two tuples coincide, proving uniqueness of the finite-rank solution tuple for the fixed fraction pair $Q$.

\end{proof}

We remark that even for bounded rational elements, the existence of a finite-rank solution to equation~\eqref{eq:certificate} can depend on the fraction pair; see Example~\ref{ex:choice-fraction-pair}.
We end this section by a remark on the connection between our criteria and Hochschild homology.

\begin{remark}\label{rem:boundary-interpretation}
Equip $L^0(\M)\otimes_{\C}L^0(\M)^{\mathrm{op}}$ with the $L^0(\M)$-bimodule structure defined, for $u,v\in L^0(\M)$ and $k\in L^0(\M)\otimes_{\C}L^0(\M)^{\mathrm{op}}$, by
\[
 u\cdot k\cdot v:=k(v\otimes u^{\mathrm{op}}).
\]
The degree-one Hochschild differential with these coefficients \cite[Section~9.1.1, p.~300]{Wei94} is the map
\[
 d_1:\bigl(L^0(\M)\otimes_{\C}L^0(\M)^{\mathrm{op}}\bigr)
 \otimes_{\C}L^0(\M)
 \longrightarrow L^0(\M)\otimes_{\C}L^0(\M)^{\mathrm{op}}
\]
given by
\[
 d_1(k\otimes u):=k\cdot u-u\cdot k
 =k(u\otimes1-1\otimes u^{\mathrm{op}}).
\]
Put $V_X:=\operatorname{span}_{\C}\{X_1,\ldots,X_d\}$.
Restricting the last tensor factor to $V_X$, we obtain
\[
 \partial_X:
 \bigl(L^0(\M)\otimes_{\C}L^0(\M)^{\mathrm{op}}\bigr)
 \otimes_{\C}V_X
 \longrightarrow L^0(\M)\otimes_{\C}L^0(\M)^{\mathrm{op}},
 \qquad \partial_X:=d_1|_{(L^0(\M)\otimes_{\C}L^0(\M)^{\mathrm{op}})\otimes_{\C}V_X}.
\]
By Definition~\ref{def:psi},
\[
 \partial_X\!\left(\sum_i k_i\otimes X_i\right)
 =\sum_i\Psi_i(k_i).
\]
This restriction to the linear span of the generators follows \cite[Section~3.1.1, equation~(3.1)]{CS05}, with $L^0(\M)\otimes_{\C}L^0(\M)^{\mathrm{op}}$ as the coefficient bimodule.
Thus condition~\textup{(\ref{item:delta-tensors})} in Theorem~\ref{thm:delta} requires $u\otimes1-1\otimes u^{\mathrm{op}}$ to lie in the image of $\partial_X$.

Similarly, equip $\F(H)$ with the $\A_X$-bimodule structure
\[
 x\cdot T\cdot y:=\rho(y)T\rho(x),
 \qquad x,y\in\A_X,\quad T\in\F(H).
\]
The degree-one Hochschild differential is then
\[
 d_1:\F(H)\otimes_{\C}\A_X\longrightarrow\F(H),
 \qquad d_1(T\otimes x):=T\cdot x-x\cdot T=[\rho(x),T].
\]
Since $\A_X$ is generated by $X_1,\ldots,X_d$, the identity
\[
 [\rho(xy),T]
 =[\rho(y),\rho(x)T]+[\rho(x),T\rho(y)]
\]
shows, by induction on monomials, that restricting the last tensor factor to $V_X$ gives the full image of $d_1$:
\[
 \operatorname{im}d_1
 =d_1\bigl(\F(H)\otimes_{\C}V_X\bigr)
 =\left\{\sum_{i=1}^d[\rho(X_i),T_i]:T_i\in\F(H)\right\}.
\]
Condition \eqref{eq:certificate} is therefore
\[
 d_1\!\left(\sum_iT_i\otimes X_i\right)=\Gamma[Q,P_\Omega].
\]
Thus $\Gamma[Q,P_\Omega]$ is a boundary, or equivalently its class vanishes in $H_0(\A_X,\F(H)):=\F(H)/\operatorname{im}d_1$.
\end{remark}

\section{Liberation dual system criterion}
\label{sec:liberation-dual-system-criterion}

We prove Theorem~\ref{thm:intrinsic-analytic} in this section.
We begin by recalling Definition~\ref{def:liberation-dual-system}.

\begin{definition*}
Let $X=(X_1,\ldots,X_d)\in\M^d$, and put $H:=L^2(\M,\tau)$ and $\Omega:=1$.
A \emph{liberation dual system} for $X$ is a tuple $(L_1,\ldots,L_d)$ of linear operators on $H$ such that
\[
 \M\Omega\subseteq\dom L_i\cap\dom L_i^*
 \qquad(1\leq i\leq d)
\]
and, on $\M\Omega$,
\[
 [X_j,L_i]=\delta_{ij}[X_i,P_\Omega]
 \qquad(1\leq i,j\leq d).
\]
Here $\M\Omega:=\{x\Omega:x\in\M\}$ and $P_\Omega$ is the orthogonal projection onto $\C\Omega$.
We call this system \emph{normalized} if, on $\M\Omega$,
\[
 \sum_{i=1}^dL_i=P_\Omega.
\]
\end{definition*}

\begin{remark}\label{rem:liberation-gradients}
For $X\in\M^d$, put $A_i:=W^*(X_i)$ and $B_i:=W^*(X_j:j\ne i)$, with $B_1:=\C1$ if $d=1$.
If the liberation gradients $j(A_i:B_i)$ in the sense of \cite[Sections~5.3--5.4]{Voi99} exist in $L^2(\M,\tau)$ for every $i$, then $X$ admits a normalized liberation dual system in the sense of Definition~\ref{def:liberation-dual-system}.
\end{remark}

Let $L$ be a linear operator on $H$ with $\M\Omega\subseteq\dom L$.
For $c,e\in\M$, define on $\M\Omega$
\[
 m_L(c\otimes e^{\mathrm{op}}):=cLe.
\]
Since $e\M\Omega\subseteq\M\Omega$, this expression is well defined and bilinear in $c$ and $e$.
Extending linearly gives $m_L(\ell):\M\Omega\to H$ for every $\ell\in\M\otimes_{\C}\M^{\mathrm{op}}$.

Let $(L_1,\ldots,L_d)$ be a liberation dual system for $X$.
Then for $\ell\in\M\otimes_{\C}\M^{\mathrm{op}}$ and $1\leq i,j\leq d$, we have on $\M\Omega$
\begin{equation}\label{eq:liberation-tensor-actions}
 m_{L_i}(\Psi_j(\ell))=\delta_{ij}\Theta(\Psi_i(\ell)).
\end{equation}
Indeed, for $\ell=c\otimes d^{\mathrm{op}}$ with $c,d\in\M$, we have, on $\M\Omega$,
\[
 \begin{aligned}
 m_{L_i}(\Psi_j(\ell))
 &=m_{L_i}(cX_j\otimes d^{\mathrm{op}}-c\otimes(X_jd)^{\mathrm{op}})\\
 &=cX_jL_id-cL_iX_jd\\
 &=c[X_j,L_i]d\\
 &=\delta_{ij}c[X_i,P_\Omega]d\\
 &=\delta_{ij}\Theta(\Psi_i(\ell)).
 \end{aligned}
\]
The fourth equality uses Definition~\ref{def:liberation-dual-system}, and the last uses equation~\eqref{eq:theta-bounded}.
Linearity gives equation~\eqref{eq:liberation-tensor-actions} for every $\ell$.

For $\ell\in\M\otimes_{\C}\M^{\mathrm{op}}$ and $s,t\in\M$, the definition of $m_{L_i}$ also gives on $\M\Omega$
\begin{equation}\label{eq:liberation-phi-bimodule}
 m_{L_i}((s\otimes t^{\mathrm{op}})\ell)=sm_{L_i}(\ell)t.
\end{equation}
For $\ell=c\otimes d^{\mathrm{op}}$, multiplication in the opposite factor gives
\[
 (s\otimes t^{\mathrm{op}})(c\otimes d^{\mathrm{op}})
 =sc\otimes(dt)^{\mathrm{op}},
\]
and hence, on $\M\Omega$,
\[
 m_{L_i}((s\otimes t^{\mathrm{op}})(c\otimes d^{\mathrm{op}}))
 =scL_idt
 =sm_{L_i}(c\otimes d^{\mathrm{op}})t.
\]
Linearity proves equation~\eqref{eq:liberation-phi-bimodule} for every $\ell$.

\begin{proposition}\label{prop:liberation-normalization}
If $X\in\M^d$ admits a liberation dual system, then it admits a normalized liberation dual system.
\end{proposition}
\begin{proof}
Let $(L_1,\ldots,L_d)$ be a liberation dual system for $X$.
Define operators $C,\widetilde L_1,\ldots,\widetilde L_d$ on $\M\Omega$ by
\[
 C:=\sum_{k=1}^dL_k-P_\Omega,
 \qquad \widetilde L_i:=L_i-\frac1dC.
\]
Since $X_j\M\Omega\subseteq\M\Omega$, the commutator $[X_j,C]$ is defined on $\M\Omega$.
On this subspace,
\[
 [X_j,C]
 =\sum_{k=1}^d[X_j,L_k]-[X_j,P_\Omega]
 =\sum_{k=1}^d\delta_{kj}[X_k,P_\Omega]-[X_j,P_\Omega]
 =0.
\]
The second equality uses Definition~\ref{def:liberation-dual-system}.
Consequently, for $\xi\in\M\Omega$ and $1\leq i,j\leq d$,
\[
 [X_j,\widetilde L_i]\xi
 =[X_j,L_i]\xi=\delta_{ij}[X_i,P_\Omega]\xi,
 \qquad
 \sum_{i=1}^d\widetilde L_i\xi
 =\sum_{i=1}^dL_i\xi-C\xi
 =P_\Omega\xi.
\]
Since $\mathcal M\Omega\subseteq\operatorname{dom}L_i^*$, we have for $\eta,\xi\in\M\Omega$,
\[
 \langle\eta,\widetilde L_i\xi\rangle
 =\left\langle
 L_i^*\eta-\frac1d\sum_{k=1}^dL_k^*\eta+\frac1dP_\Omega\eta,
 \xi\right\rangle.
\]
By the definition of the adjoint, \(\eta\in\operatorname{dom}\widetilde L_i^*\), and hence \(\mathcal M\Omega\subseteq\operatorname{dom}\widetilde L_i^*\).
Thus \(\widetilde L_1,\ldots,\widetilde L_d\) satisfy all the defining conditions.
\end{proof}

\begin{lemma}[{see \cite[Problem~III.5.18]{Kat95}}]\label{lem:finite-closure}
Let \(T\) be a densely defined closable operator on a Hilbert space \(H\).
If \(T\) has finite rank, then \(\overline T\in\F(H)\), with the same range and rank as \(T\).
\end{lemma}

\begin{lemma}\label{lem:tensor-difference-injectivity}
Let $V,W$ be complex vector spaces, and let $A:V\to V$ and $B:W\to W$ be linear maps.
If $A$ or $B$ has no nonzero finite-dimensional invariant subspace, then
\[
 A\otimes I_W-I_V\otimes B:V\otimes_{\C}W\longrightarrow V\otimes_{\C}W
\]
is injective.
\end{lemma}
\begin{proof}
Suppose that $B$ has no nonzero finite-dimensional invariant subspace.
If a nonzero kernel vector existed, write it as $\sum_{r=1}^n v_r\otimes w_r$ with $n$ minimal.
By minimality of $n$, both families are linearly independent.
Set $E:=\operatorname{span}\{w_1,\ldots,w_n\}$, and let $\pi:W\to W/E$ be the quotient map.
Since the tensor lies in the kernel, we have
\[
 \sum_{r=1}^n Av_r\otimes w_r=\sum_{r=1}^n v_r\otimes Bw_r.
\]
Applying $I_V\otimes \pi$ and using $\pi(w_r)=0$ gives
\[
 \sum_{r=1}^n v_r\otimes \pi(Bw_r)=0.
\]
Since $v_1,\ldots,v_n$ are linearly independent, choose linear functionals $\varphi_j:V\to\C$ satisfying $\varphi_j(v_r)=\delta_{jr}$.
Applying $\varphi_j\otimes I_{W/E}$ yields $\pi(Bw_j)=0$, and hence $Bw_j\in E$, for every $1\leq j\leq n$.
Thus $B(E)\subseteq E$, contradicting the assumption on $B$, since $E$ is nonzero and finite-dimensional.
The case where $A$ has the stated property follows by interchanging the tensor factors.
\end{proof}

\begin{lemma}\label{lem:tensor-quotient-kernels}
Let $E$ be a subspace of a complex vector space $V$, and let $\pi:V\to V/E$ be the quotient map.
For every complex vector space $W$,
\[
 \ker(\pi\otimes I_W)=E\otimes_{\C}W,
 \qquad
 \ker(I_W\otimes\pi)=W\otimes_{\C}E.
\]
\end{lemma}
\begin{proof}
Choose a complement $V_1$ such that $V=E\oplus V_1$.
Then
\[
 V\otimes_{\C}W
 =(E\otimes_{\C}W)\oplus(V_1\otimes_{\C}W),
\]
and the map $\pi\otimes I_W$ vanishes on the first summand and maps the second summand isomorphically onto $(V/E)\otimes_{\C}W$, since $\pi|_{V_1}:V_1\to V/E$ is an isomorphism.
This proves the first kernel equality.
Interchanging the tensor factors proves injectivity of the other natural map and the second kernel equality.
\end{proof}

\begin{lemma}\label{lem:tensor-l2-recovery}
Let $x\in\M$ and suppose that $x-\lambda$ is regular for every $\lambda\in\C$.
The map $\Psi_x$ of Definition~\ref{def:psi} is injective and satisfies
\[
 \Psi_x^{-1}(H\otimes_{\C}H)=H\otimes_{\C}H.
\]
\end{lemma}
\begin{proof}
Put $U:=L^0(\M)$ and define the right multiplication maps
\[
 \begin{aligned}
 R_x:U&\longrightarrow U,& R_x(u)&:=ux,\\
 R_{x^{\mathrm{op}}}:U^{\mathrm{op}}&\longrightarrow U^{\mathrm{op}},&
 R_{x^{\mathrm{op}}}(v^{\mathrm{op}})&:=(xv)^{\mathrm{op}}.
 \end{aligned}
\]
Thus $\Psi_x=R_x\otimes I-I\otimes R_{x^{\mathrm{op}}}$.
Neither right multiplication map has a nonzero finite-dimensional invariant subspace: such a subspace would contain an eigenvector, contradicting regularity of $x-\lambda$.
Lemma~\ref{lem:tensor-difference-injectivity} therefore implies that $\Psi_x$ is injective.

For the preimage equality, observe that $x\in\M$ implies that both maps preserve $H$.
Their restrictions to $H$ also have no nonzero finite-dimensional invariant subspace.
In particular, $\Psi_x$ preserves $H\otimes_{\C}H$.

Let $\pi_1:U\to U/H$ and $\pi_2:U^{\mathrm{op}}\to U^{\mathrm{op}}/H$ be the quotient maps.
The induced linear maps on the quotient spaces are
\[
 \begin{aligned}
 \widetilde R_x:U/H&\longrightarrow U/H,&
 \widetilde R_x(u+H)&:=ux+H,\\
 \widetilde R_{x^{\mathrm{op}}}:U^{\mathrm{op}}/H&\longrightarrow U^{\mathrm{op}}/H,&
 \widetilde R_{x^{\mathrm{op}}}(v^{\mathrm{op}}+H)&:=(xv)^{\mathrm{op}}+H.
 \end{aligned}
\]
They satisfy
\[
 \pi_1R_x=\widetilde R_x\pi_1,
 \qquad
 \pi_2R_{x^{\mathrm{op}}}=\widetilde R_{x^{\mathrm{op}}}\pi_2.
\]
Now let $k\in U\otimes_{\C}U^{\mathrm{op}}$ satisfy $\Psi_x(k)\in H\otimes_{\C}H$.
Using $\pi_1R_x=\widetilde R_x\pi_1$, we obtain
\[
 \begin{aligned}
 &\bigl(\widetilde R_x\otimes I-I\otimes R_{x^{\mathrm{op}}}\bigr)
    (\pi_1\otimes I)k\\
 &\qquad=(\pi_1\otimes I)\Psi_x(k)=0.
 \end{aligned}
\]
The last equality follows from $\Psi_x(k)\in H\otimes_{\C}H$ and $\pi_1|_H=0$.
By Lemma~\ref{lem:tensor-difference-injectivity}, the map $\widetilde R_x\otimes I-I\otimes R_{x^{\mathrm{op}}}$ is injective on $(U/H)\otimes_{\C}U^{\mathrm{op}}$, since $R_{x^{\mathrm{op}}}$ has no nonzero finite-dimensional invariant subspace.
Consequently, $(\pi_1\otimes I)k=0$.
Applying the first kernel equality of Lemma~\ref{lem:tensor-quotient-kernels} with $V=U$, $E=H$, $W=U^{\mathrm{op}}$ and $\pi=\pi_1$ gives
\[
 k\in\ker(\pi_1\otimes I)=H\otimes_{\C}U^{\mathrm{op}}.
\]

Now regard $k$ as an element of $H\otimes_{\C}U^{\mathrm{op}}$.
Using $\pi_2R_{x^{\mathrm{op}}}=\widetilde R_{x^{\mathrm{op}}}\pi_2$, we similarly obtain
\[
 k\in\ker(I\otimes \pi_2)=H\otimes_{\C}H.
\]
\end{proof}

To prove Theorem~\ref{thm:intrinsic-analytic} using Theorem~\ref{thm:delta}, we first need to establish that the existence of a liberation dual system implies the maximality of $\Delta$, which is captured by the following proposition.

\begin{proposition}\label{prop:liberation-delta}
Let $X=(X_1,\ldots,X_d)\in\M^d$, and suppose that $X_i-\lambda$ is regular for every $1\leq i\leq d$ and every $\lambda\in\C$.
If $X$ admits a liberation dual system, then $\Delta(X)=d$.
\end{proposition}
\begin{proof}
Let $T_1,\ldots,T_d\in\F(H)$ satisfy $\sum_i[\rho(X_i),T_i]=0$, and put $k_i:=\Theta^{-1}(T_i)\in H\otimes_{\C}H$.
Since $X_i\in\M$, each $\Psi_i(k_i)$ belongs to $H\otimes_{\C}H$.
By equation~\eqref{eq:tensor-commutator-general},
\[
 \Theta\!\left(\sum_i\Psi_i(k_i)\right)
 =\sum_i\Theta(\Psi_i(k_i))
 =\sum_i[\rho(X_i),T_i]=0.
\]
Since $\Theta$ is injective, $\sum_i\Psi_i(k_i)=0$.
By Lemma~\ref{lem:tensor-clearing}, choose regular $s,t\in\M$ such that $\ell_i:=(s\otimes t^{\mathrm{op}})k_i\in\M\otimes_{\C}\M^{\mathrm{op}}$ for every $i$.
Equation~\eqref{eq:psi-left-linearity} yields $\sum_i\Psi_i(\ell_i)=0$.
Applying $m_{L_i}$ and equation~\eqref{eq:liberation-tensor-actions} gives, on $\M\Omega$,
\[
 0=m_{L_i}\!\left(\sum_j\Psi_j(\ell_j)\right)
 =\Theta(\Psi_i(\ell_i)).
\]
Since $\Theta(\Psi_i(\ell_i))$ is bounded and $\M\Omega$ is dense in $H$, it vanishes on $H$.
Hence $\Psi_i(\ell_i)=0$.
Since $X_i-\lambda$ is regular for every $\lambda\in\C$, Lemma~\ref{lem:tensor-l2-recovery} gives $\ell_i=0$.
Cancellation of $s\otimes t^{\mathrm{op}}$ gives $k_i=0$, and therefore $T_i=0$.
Remark~\ref{rem:cycle-characterization} now gives $\Delta(X)=d$.
\end{proof}

We record two technical lemmas that will be useful in the subsequent proof for Theorem~\ref{thm:intrinsic-analytic}.

\begin{lemma}\label{lem:tensor-images-l2}
Let $X=(X_1,\ldots,X_d)\in\M^d$, and suppose that $X_i-\lambda$ is regular for every $1\leq i\leq d$ and every $\lambda\in\C$.
Let $(L_1,\ldots,L_d)$ be a liberation dual system for $X$.
If $k_1,\ldots,k_d\in L^0(\M)\otimes_{\C}L^0(\M)^{\mathrm{op}}$ satisfy
\[
 \ell_0:=\sum_{i=1}^d\Psi_i(k_i)\in\M\otimes_{\C}\M^{\mathrm{op}},
\]
then, for every $i$, $k_i\in H\otimes_{\C}H$ and
\begin{equation}\label{eq:liberation-tensor-closure}
 \overline{m_{L_i}(\ell_0)}=\Theta(\Psi_i(k_i)).
\end{equation}
\end{lemma}
\begin{proof}
Fix $i$, write $\ell_0=\sum_r c_r\otimes d_r^{\mathrm{op}}$ with $c_r,d_r\in\M$, and put $T_i:=m_{L_i}(\ell_0)$ on $\M\Omega$.
By Definition~\ref{def:liberation-dual-system}, $\M\Omega\subseteq\dom L_i^*$.
Since $c_r^*\M\Omega\subseteq\M\Omega$, the definition of the adjoint gives $\M\Omega\subseteq\dom T_i^*$ and, on $\M\Omega$,
\[
 T_i^*=\sum_r d_r^*L_i^*c_r^*.
\]
Thus $T_i$ is closable.
By Lemma~\ref{lem:tensor-clearing}, choose regular $s,t\in\M$ such that
\[
 \ell_j:=(s\otimes t^{\mathrm{op}})k_j
 \in\M\otimes_{\C}\M^{\mathrm{op}}
 \qquad(1\leq j\leq d).
\]
Equation~\eqref{eq:psi-left-linearity} gives
\[
 \sum_j\Psi_j(\ell_j)=(s\otimes t^{\mathrm{op}})\ell_0.
\]
Using equation~\eqref{eq:liberation-phi-bimodule} in the second equality below and equation~\eqref{eq:liberation-tensor-actions} in the last equality, we obtain, on $\M\Omega$,
\begin{equation}\label{eq:liberation-cleared-operator}
 \begin{aligned}
 sT_it
 &=sm_{L_i}(\ell_0)t\\
 &=m_{L_i}((s\otimes t^{\mathrm{op}})\ell_0)\\
 &=m_{L_i}\!\left(\sum_j\Psi_j(\ell_j)\right)\\
 &=\Theta(\Psi_i(\ell_i)).
 \end{aligned}
\end{equation}
The operator $\Theta(\Psi_i(\ell_i))$ has finite rank and $s$ is injective, so the restriction of $T_i$ to $t\M\Omega$ has finite rank.
Since $t$ has dense range, $t\M\Omega$ is dense in $H$.
This restriction is closable, and Lemma~\ref{lem:finite-closure} makes its closure a bounded finite-rank operator on $H$.
The operator $\overline{T_i}$ extends that closure and therefore equals it.
Set $h_i:=\Theta^{-1}(\overline{T_i})\in H\otimes_{\C}H$.
Extending equation~\eqref{eq:liberation-cleared-operator} by boundedness and density, we have
\[
 s\overline{T_i}t=\Theta(\Psi_i(\ell_i)).
\]
Applying $\Theta^{-1}$ to both sides, we obtain
\[
 (s\otimes t^{\mathrm{op}})h_i
 =\Psi_i(\ell_i)
 =(s\otimes t^{\mathrm{op}})\Psi_i(k_i).
\]
Cancellation of the regular outside factors gives $\Psi_i(k_i)=h_i\in H\otimes_{\C}H$.
Since every $X_i-\lambda$ is regular, Lemma~\ref{lem:tensor-l2-recovery} gives $k_i\in H\otimes_{\C}H$.
\end{proof}

\begin{lemma}\label{lem:fraction-domains}
Let \(D_1,\ldots,D_d\) be densely defined linear operators on \(H\) such that
\[\M\Omega\subseteq\dom D_i\cap\dom D_i^*\]
for \(1\leq i\leq d\).
For every \(u\in L^0(\M)\), every fraction pair \(Q\) of \(u\), and every \(i\), the operator \(\Gamma[Q,D_i]\) is closable and densely defined.
\end{lemma}
\begin{proof}
Fix \(i\) and write \(Q=(s,a,b,t)\).
Since \(b,t\in\M\) preserve \(\M\Omega\), the domain of \(\Gamma[Q,D_i]\) contains \(\M\Omega\).
For \(\eta\in\M\Omega\), we have \(s^*\eta,a^*\eta\in\M\Omega\subseteq\dom D_i^*\).
Thus, for \(\xi\in\dom\Gamma[Q,D_i]\),
\begin{equation}\label{eq:defect-adjoint}
 \langle\eta,\Gamma[Q,D_i]\xi\rangle
 =\left\langle
 b^*D_i^*(s^*\eta)-t^*D_i^*(a^*\eta),\xi
 \right\rangle.
\end{equation}
Thus \(\M\Omega\subseteq\dom(\Gamma[Q,D_i])^*\), so \(\Gamma[Q,D_i]\) is closable.
\end{proof}

\begin{proof}[Proof of Theorem~\ref{thm:intrinsic-analytic}]
By Proposition~\ref{prop:liberation-delta}, $\Delta(X)=d$.
We prove $(1)\implies(4)\implies(3)\implies(2)\implies(1)$.

$(1)\implies(4)$.
Fix any fraction pair $Q=(s,a,b,t)$ for $u$.
By Theorem~\ref{thm:delta}\textup{(\ref{item:delta-all-fraction-tensors})}, there are tensors $\ell_i\in L^0(\M)\otimes_{\C}L^0(\M)^{\mathrm{op}}$ satisfying equation~\eqref{eq:all-fraction-tensor-criterion}.
Set
\[
 \ell_0:=a\otimes t^{\mathrm{op}}-s\otimes b^{\mathrm{op}}
   =\sum_i\Psi_i(\ell_i).
\]
Lemma~\ref{lem:tensor-images-l2} gives $\ell_i\in H\otimes_{\C}H$ and $\overline{m_{L_i}(\ell_0)}=\Theta(\Psi_i(\ell_i))$.
For $\xi\in\M\Omega$, by the definition of $m_{L_i}$ and $\ell_0$,
\[
 m_{L_i}(\ell_0)\xi
 = aL_i(t\xi)-sL_i(b\xi)
 = -\Gamma[Q,L_i]\xi.
\]
By Definition~\ref{def:liberation-dual-system}, $\M\Omega\subseteq\dom L_i\cap\dom L_i^*$.
Lemma~\ref{lem:fraction-domains} therefore implies that $\Gamma[Q,L_i]$ is closable.
Thus we have
\[
\overline{\Gamma[Q,L_i]}=-\overline{m_{L_i}(\ell_0)}=-\Theta(\Psi_i(\ell_i)).
\]
With $T_i:=-\Theta(\ell_i)\in\F(H)$, equation~\eqref{eq:tensor-commutator-general} now gives
\[
 \overline{\Gamma[Q,L_i]}
 =-\Theta(\Psi_i(\ell_i))
 =[\rho(X_i),T_i].
\]
This gives equation~\eqref{eq:intrinsic-analytic-coordinates} and proves assertion~\textup{(4)}.

$(4)\implies(3)$.
Fix any fraction pair $Q=(s,a,b,t)$ for $u$ and choose $T_1,\ldots,T_d\in\F(H)$ satisfying equation~\eqref{eq:intrinsic-analytic-coordinates}.
By equation~\eqref{eq:liberation-normalization}, on $\M\Omega$ we have
\[
 \begin{aligned}
  \sum_i[\rho(X_i),T_i]
  &=\sum_i\overline{\Gamma[Q,L_i]}\\
  &=s\Bigl(\sum_iL_i\Bigr)b-a\Bigl(\sum_iL_i\Bigr)t\\
  &=sP_\Omega b-aP_\Omega t
   =\Gamma[Q,P_\Omega].
 \end{aligned}
\]
The bounded operators $\sum_i[\rho(X_i),T_i]$ and $\Gamma[Q,P_\Omega]$ agree on the dense subspace $\M\Omega$ and hence on all of $H$.
This proves assertion~\textup{(3)}.

$(3)\implies(2)$.
Let $Q$ be the fraction pair for $u$ given by equation~\eqref{eq:canonical}.
By assertion~\textup{(3)}, there exist $T_1,\ldots,T_d\in\F(H)$ satisfying equation~\eqref{eq:intrinsic-analytic-criterion} for this $Q$.

$(2)\implies(1)$.
This is the implication $\textup{(\ref{item:delta-finite-rank})}\implies\textup{(1)}$ of Theorem~\ref{thm:delta}.

Theorem~\ref{thm:delta} gives uniqueness of the solution of equation~\eqref{eq:intrinsic-analytic-criterion} for a fixed fraction pair.
The implication $(4)\implies(3)$ gives uniqueness for equation~\eqref{eq:intrinsic-analytic-coordinates} as well.
\end{proof}

\section{Dual system criterion}
\label{sec:dual-system-criterion}

We prove Theorem~\ref{thm:analytic} by applying Theorem~\ref{thm:intrinsic-analytic} to a normalized liberation dual system constructed from the given dual system.
We first recall Definition~\ref{def:dual-system} of the dual system for reader's convenience.

\begin{definition*}
Let $X=(X_1,\ldots,X_d)\in\M^d$.
A \emph{dual system} for $X$ is a tuple $(D_1,\ldots,D_d)$ of linear operators on $H$ such that
\[
 \M\Omega\subseteq\dom D_i\cap\dom D_i^*,\qquad
 [X_j,D_i]|_{\M\Omega}=\delta_{ij}P_\Omega|_{\M\Omega}
 \qquad(1\leq i,j\leq d).
\]
\end{definition*}

A selfadjoint tuple generating $\M$ admits a dual system in our sense if and only if it admits one in the sense of Shlyakhtenko \cite[Theorem~1]{Shl06}.
For the reverse implication, one may use the estimates in \cite[Lemma~3.4]{MSW17} and Kaplansky density theorem to obtain both domain inclusions in our definition.

\begin{remark}\label{rem:dual-maximality}
If $X$ admits a dual system in the sense of Definition~\ref{def:dual-system}, then $\Delta(X)=d$.
The argument in \cite[proof of Proposition~6.2]{MSY23} applies \emph{mutatis mutandis}.
\end{remark}

We first prove that the existence of a finite-rank closure and its rank are independent of the fraction pair, and identify finite-rank solutions of equation~\eqref{eq:certificate} with the closures of the corresponding fraction commutators.
For the right annihilation operators of a free semicircular system, the finite-rank equivalence in Lemma~\ref{lem:fraction-ranks} is the one proved by Miyagawa \cite[Lemma~4.6]{Miy22}.

\begin{lemma}\label{lem:fraction-ranks}
Let $D$ be a densely defined linear operator on $H$ such that
$\M\Omega\subseteq\dom D\cap\dom D^*$.
For $u\in L^0(\M)$ and any two fraction pairs $Q,Q'$ for $u$,
\[
 \overline{\Gamma[Q,D]}\in\F(H)
 \quad\Longleftrightarrow\quad
 \overline{\Gamma[Q',D]}\in\F(H).
\]
When these conditions hold, the two closures have the same rank.
\end{lemma}
\begin{proof}
Write $Q=(s,a,b,t)$ and $Q'=(s',a',b',t')$.
Lemma~\ref{lem:common} gives regular $x,y,p,q\in\M$ such that $xs=ys'$ and $tp=t'q$.
Since both fraction pairs represent $u$, we have $a=su$, $a'=s'u$, $b=ut$ and $b'=ut'$ in $L^0(\M)$.
Thus
\begin{align*}
 xa&=x(su)=(xs)u=(ys')u=y(s'u)=ya',\\
 bp&=(ut)p=u(tp)=u(t'q)=(ut')q=b'q.
\end{align*}
By equation~\eqref{eq:fraction-commutator}, these identities and the common multiples give on $\M\Omega$
\begin{equation}\label{eq:fraction-common-comparison}
 \begin{aligned}
 x\Gamma[Q,D]p
 &=(xs)D(bp)-(xa)D(tp)
  =(ys')D(b'q)-(ya')D(t'q)\\
 &=y\bigl(s'D(b'q)-a'D(t'q)\bigr)
  =y\Gamma[Q',D]q.
 \end{aligned}
\end{equation}
Suppose that $\overline{\Gamma[Q,D]}\in\F(H)$.
The left-hand side of equation~\eqref{eq:fraction-common-comparison} has finite rank, and $y$ is injective.
Thus the restriction of $\Gamma[Q',D]$ to the dense subspace $q\M\Omega$ has finite rank.
It is closable by Lemma~\ref{lem:fraction-domains}, so Lemma~\ref{lem:finite-closure} makes its closure a bounded finite-rank operator on all of $H$.
The closure of $\Gamma[Q',D]$ extends this operator and hence equals it.
Interchanging $Q$ and $Q'$ proves the equivalence.

Both closures are now bounded, so equation~\eqref{eq:fraction-common-comparison} extends to
\[
 x\overline{\Gamma[Q,D]}p=y\overline{\Gamma[Q',D]}q
 \qquad\text{on }H.
\]
The left factors $x,y$ are injective, while $p,q$ have dense range.
These multiplications preserve the rank of a bounded finite-rank operator, proving the rank equality.
\end{proof}

We next show that, for a dual system $(D_1,\ldots,D_d)$, every finite-rank solution $(T_1,\ldots,T_d)$ of equation~\eqref{eq:certificate} satisfies $T_i=\overline{\Gamma[Q,D_i]}$ for every $i$.

\begin{lemma}\label{lem:finite-rank-identification}
Let $X=(X_1,\ldots,X_d)\in\M^d$, and let $(D_1,\ldots,D_d)$ be a dual system for $X$.
Let $Q$ be a fraction pair for $u\in L^0(\M)$.
If $T_1,\ldots,T_d\in\F(H)$ satisfy
\[
 \sum_j[\rho(X_j),T_j]=\Gamma[Q,P_\Omega],
\]
then $\overline{\Gamma[Q,D_i]}=T_i$ for every $i$.
\end{lemma}
\begin{proof}
Write $Q=(s,a,b,t)$ and put $\ell_j:=-\Theta^{-1}(T_j)\in H\otimes_{\C}H$.
By Lemma~\ref{lem:tensor-certificate},
\[
 \sum_j\Psi_j(\ell_j)=a\otimes t^{\mathrm{op}}-s\otimes b^{\mathrm{op}}.
\]
Lemma~\ref{lem:tensor-clearing} gives regular $s',t'\in\M$ such that
$\ell_j':=(s'\otimes(t')^{\mathrm{op}})\ell_j\in\M\otimes_{\C}\M^{\mathrm{op}}$ for every $j$.
Equation~\eqref{eq:psi-left-linearity} gives
\begin{equation}\label{eq:dual-certificate-tensors}
 \sum_j\Psi_j(\ell_j')
 =s'a\otimes(tt')^{\mathrm{op}}-s's\otimes(bt')^{\mathrm{op}}.
\end{equation}
Fix $i$ and $\xi\in\M\Omega$, and apply $\ell\mapsto m_{D_i}(\ell)\xi$ to both sides of equation~\eqref{eq:dual-certificate-tensors}.
To evaluate the left-hand side of equation~\eqref{eq:dual-certificate-tensors}, we first calculate on a simple tensor $c\otimes e^{\mathrm{op}}$ with $c,e\in\M$.
Expanding $\Psi_j$ gives
\[
 \begin{aligned}
 m_{D_i}(\Psi_j(c\otimes e^{\mathrm{op}}))\xi
 &=cX_jD_i(e\xi)-cD_i(X_je\xi)\\
 &=\delta_{ij}cP_\Omega e\xi
  =\delta_{ij}\Theta(c\otimes e^{\mathrm{op}})\xi.
 \end{aligned}
\]
The second equality uses Definition~\ref{def:dual-system} of a dual system at $e\xi\in\M\Omega$, and the last uses equation~\eqref{eq:theta-bounded}.
By linearity, we obtain
\begin{align*}
 \sum_jm_{D_i}(\Psi_j(\ell_j'))\xi
 &=\sum_j\delta_{ij}\Theta(\ell_j')\xi\\
 &=\Theta(\ell_i')\xi\\
 &=s'\Theta(\ell_i)(t')\xi\\
 &=-s'T_it'\xi.
\end{align*}
Evaluating the right-hand side of equation~\eqref{eq:dual-certificate-tensors}, we obtain
\[
s'aD_i(tt'\xi)-s'sD_i(bt'\xi)=-s'\Gamma[Q,D_i]t'\xi.
\]
Since $s'$ is injective, $\Gamma[Q,D_i]$ agrees with $T_i$ on the dense subspace $t'\M\Omega$.
Lemma~\ref{lem:fraction-domains} gives closability of $\Gamma[Q,D_i]$.
Its closure therefore extends the bounded operator $T_i$, whose domain is all of $H$, and hence equals $T_i$.
\end{proof}

For \(N,S\in\B(H)\) with \(N\) normal, the classical Fuglede theorem \cite[Theorem~I, p.~35]{Fug50} gives \(SN=NS\implies SN^*=N^*S\).
The following analogue for bounded linear maps from \(\M\) to \(H\) will be used in Lemma~\ref{lem:module-commutation}.

\begin{lemma}[{see \cite[Proposition~1]{CS77}}]\label{lem:normal-intertwining}
Let \(x\in\M\) be normal, and let \(T:\M\to H\) be a bounded linear map.
If \(T(xz)=xT(z)\) for every \(z\in\M\), then \(T(x^*z)=x^*T(z)\) for every \(z\in\M\).
\end{lemma}

\begin{lemma}\label{lem:module-commutation}
Let \(T:\M\Omega\to H\) be a closable linear operator.
Then \(y\mapsto T(y\Omega)\) is bounded from \(\M\) to \(H\).
Suppose \(x_1,\ldots,x_m\in\M\) are normal operators and \(Tx_j=x_jT\) on \(\M\Omega\) for every \(j\).
For every \(a\in W^*(x_1,\ldots,x_m)\), we then have on \(\M\Omega\)
\[
 Ta=aT.
\]
\end{lemma}
\begin{proof}
Suppose \(y_n\to y\) in the operator norm of \(\M\) and \(T(y_n\Omega)\to\eta\) in \(H\).
Then
\[
 \|(y_n-y)\Omega\|_2\leq\|y_n-y\|\longrightarrow 0.
\]
Since \(\overline T\) is closed and agrees with \(T\) on \(\M\Omega\), we have \(T(y\Omega)=\eta\).
Thus the map \(y\mapsto T(y\Omega)\) from \(\M\) to \(H\) has closed graph and is bounded by the closed graph theorem.
Lemma~\ref{lem:normal-intertwining} gives \(T(x_j^*\xi)=x_j^*T\xi\).
Products and linear combinations give the same relation for the unital $*$-algebra generated by \(x_1,\ldots,x_m\).
For \(a\in W^*(x_1,\ldots,x_m)\), Kaplansky density \cite[Theorem~2.3.1]{AP18} gives a uniformly bounded net \(a_\lambda\) in this algebra converging strongly to \(a\).
For fixed \(\xi\in\M\Omega\),
\[
 a_\lambda\xi\longrightarrow a\xi,\qquad
 T(a_\lambda\xi)=a_\lambda T\xi\longrightarrow aT\xi.
\]
Closedness of \(\overline T\) and \(a\xi\in\M\Omega\) give \(T(a\xi)=aT\xi\).
\end{proof}


We adapt Voiculescu's construction in \cite[Remark~13.2(c)]{Voi99} to obtain a normalized liberation dual system, whenever a dual system exists.

\begin{proposition}\label{prop:dual-liberation-construction}
Let $X\in\M^d$ and let $(D_1,\ldots,D_d)$ be a dual system for $X$.
Define on $\M\Omega$
\begin{equation}\label{eq:dual-liberation-construction}
 K_i:=[\rho(X_i),D_i|_{\M\Omega}],\qquad
 C:=\sum_jK_j-P_\Omega,\qquad
 L_i:=K_i-\frac1d C.
\end{equation}
The tuple $(L_1,\ldots,L_d)$ is a normalized liberation dual system for $X$.
The operator $C$ is closable and satisfies $[X_j,C]=0$ on $\M\Omega$ for every $j$.
\end{proposition}
\begin{proof}
By Definition~\ref{def:dual-system}, $\M\Omega\subseteq\dom D_i^*$ for every $i$.
Since $\rho(X_i)^*=\rho(X_i^*)$ preserves $\M\Omega$, for $\xi,\eta\in\M\Omega$ we have
\[
 \langle K_i\xi,\eta\rangle
 =\big\langle \xi,D_i^*\rho(X_i^*)\eta-\rho(X_i^*)D_i^*\eta\big\rangle.
\]
Thus $\M\Omega\subseteq\dom K_i^*$ for every $i$.
By equation~\eqref{eq:dual-liberation-construction}, $C$ and each $L_i$ are linear combinations of the $K_j$ and the bounded projection $P_\Omega$.
Hence $\M\Omega\subseteq\dom C^*\cap\dom L_i^*$ for every $i$.
In particular, $C$ is closable.

Next, we want to show that $[X_j,C]=0$ on $\M\Omega$ for any $j$.
Thanks to Jacobi's identity, we have on $\M\Omega$
\[
 \begin{aligned}
  [X_j,K_i]
  &=[X_j,[\rho(X_i),D_i]]\\
  &=[[X_j,\rho(X_i)],D_i]+[\rho(X_i),[X_j,D_i]]\\
  &=[\rho(X_i),[X_j,D_i]],
 \end{aligned}
\]
where the last equality is due to $[X_j,\rho(X_i)]=0$, as the left multiplication by $X_j$ and $\rho(X_i)$ commute and preserve $\M\Omega$.
We then have on $\M\Omega$
\[
 \begin{aligned}
  [X_j,K_i]
  &=\rho(X_i)[X_j,D_i]-[X_j,D_i]\rho(X_i)\\
  &=\delta_{ij}[\rho(X_i),P_\Omega].
 \end{aligned}
\]

By equation~\eqref{eq:dual-liberation-construction}, on $\M\Omega$ we have
\[
 \begin{aligned}
  [X_j,C]
  &=\sum_i[X_j,K_i]-[X_j,P_\Omega]\\
  &=\sum_i\delta_{ij}[\rho(X_i),P_\Omega]-[X_j,P_\Omega]\\
  &=[\rho(X_j),P_\Omega]-[X_j,P_\Omega]=0.
 \end{aligned}
\]
The last equality $[\rho(X_j),P_\Omega]=[X_j,P_\Omega]$ can be concluded from the identities
\[
 \rho(X_j)P_\Omega(a\Omega)=\tau(a)X_j\Omega=X_jP_\Omega(a\Omega)
\]
and
\[
 P_\Omega\rho(X_j)(a\Omega)=\tau(aX_j)\Omega
 =\tau(X_j a)\Omega=P_\Omega X_j(a\Omega).
\]

It follows that $[X_j,L_i]=[X_j,K_i]=\delta_{ij}[X_i,P_\Omega]$ and $\sum_iL_i=P_\Omega$ on $\M\Omega$.
Thus $(L_1,\ldots,L_d)$ is a normalized liberation dual system.
\end{proof}

For normal generators we can now compare the fraction commutators of the two systems.

\begin{lemma}\label{lem:fundamental}
Let $X_1,\ldots,X_d\in\M$ be normal operators generating $\M$, and let $(D_1,\ldots,D_d)$ be a dual system for $X$.
Let $(L_1,\ldots,L_d)$ be any normalized liberation dual system for $X$.
For every fraction pair $Q$ and every $\xi\in\M\Omega$, we have
\begin{equation}\label{eq:dual-liberation-fractions}
 \overline{\Gamma[Q,L_i]}\xi
 =[\rho(X_i),\overline{\Gamma[Q,D_i]}]\xi
 \qquad(1\leq i\leq d).
\end{equation}
Consequently,
\begin{equation}\label{eq:fundamental}
 \sum_{i=1}^d[\rho(X_i),\overline{\Gamma[Q,D_i]}]\xi=\Gamma[Q,P_\Omega]\xi.
\end{equation}
\end{lemma}
\begin{proof}
Let $K_i$ be as in equation~\eqref{eq:dual-liberation-construction}, and set $C_i:=K_i-L_i$ on $\M\Omega$.
The proof of Proposition~\ref{prop:dual-liberation-construction} gives $\M\Omega\subseteq\dom K_i^*$ and $[X_j,K_i]=\delta_{ij}[X_i,P_\Omega]$ on $\M\Omega$.
Since $\M\Omega\subseteq\dom L_i^*$, we have $\M\Omega\subseteq\dom C_i^*$, so $C_i$ is closable.
By Definition~\ref{def:liberation-dual-system}, we have on $\M\Omega$
\[
 [X_j,C_i]=[X_j,K_i]-[X_j,L_i]=0
 \qquad(1\leq j\leq d).
\]
Since the normal operators $X_j$ generate $\M$, Lemma~\ref{lem:module-commutation} gives $C_i(y\xi)=yC_i\xi$ for all $y\in\M$ and $\xi\in\M\Omega$.
For a fraction pair $Q=(s,a,b,t)$ for $u$, we have $sb=s(ut)=(su)t=at$ in $L^0(\M)$.
Thus on $\M\Omega$,
\[
 \Gamma[Q,C_i]=sC_i b-aC_i t=(sb-at)C_i=0.
\]
Substituting $L_i=K_i-C_i$ into equation~\eqref{eq:fraction-commutator} gives on $\M\Omega$
\begin{equation}\label{eq:dual-liberation-on-domain}
 \begin{aligned}
  \Gamma[Q,L_i]
  &=\Gamma[Q,K_i]-\Gamma[Q,C_i]
   =\Gamma[Q,K_i]\\
  &=s[\rho(X_i),D_i]b-a[\rho(X_i),D_i]t\\
  &=\rho(X_i)(sD_i b-aD_i t)-(sD_i b-aD_i t)\rho(X_i)\\
  &=\rho(X_i)\Gamma[Q,D_i]-\Gamma[Q,D_i]\rho(X_i).
 \end{aligned}
\end{equation}
By Lemma~\ref{lem:fraction-domains}, $\Gamma[Q,D_i]$ and $\Gamma[Q,L_i]$ are closable, and their closures agree with them on $\M\Omega$.
Since $\rho(X_i)$ preserves $\M\Omega$, equation~\eqref{eq:dual-liberation-on-domain} gives equation~\eqref{eq:dual-liberation-fractions}.

For $\xi\in\M\Omega$, both $b\xi$ and $t\xi$ belong to $\M\Omega$, so the identity $\sum_iL_i=P_\Omega$ on $\M\Omega$ gives
\[
 \sum_i[\rho(X_i),\overline{\Gamma[Q,D_i]}]\xi
 =\sum_i\overline{\Gamma[Q,L_i]}\xi=\Gamma[Q,P_\Omega]\xi.
\]
\end{proof}

\begin{proof}[Proof of Theorem~\ref{thm:analytic}]
Put $\M:=W^*(X)$ and use the restriction of $\tau$ to $\M$.
Remark~\ref{rem:dual-maximality} gives $\Delta(X)=d$.
By the free-field embedding in Lemma~\ref{lem:free-field}, the image $X_i-\lambda$ of the nonzero element $x_i-\lambda$ is regular for every $i$ and every $\lambda\in\C$.
Proposition~\ref{prop:dual-liberation-construction} supplies a normalized liberation dual system for $X$.
Thus Theorem~\ref{thm:intrinsic-analytic} applies.
Fix $u\in L^0(\M)$ and a fraction pair $Q$ for $u$.

$\textup{(1)}\implies\textup{(3)}$.
Theorem~\ref{thm:intrinsic-analytic}\textup{(\ref{item:intrinsic-analytic-criterion})} gives $T_1,\ldots,T_d\in\F(H)$ satisfying
\[
 \sum_i[\rho(X_i),T_i]=\Gamma[Q,P_\Omega].
\]
Lemma~\ref{lem:finite-rank-identification} gives $\overline{\Gamma[Q,D_i]}=T_i$ for every $i$, proving assertion~\textup{(3)}.

$\textup{(3)}\implies\textup{(2)}$.
A simple restriction gives \textup{(2)}.

$\textup{(2)}\implies\textup{(1)}$.
Lemmas~\ref{lem:fraction-domains} and~\ref{lem:finite-closure} give $T_i:=\overline{\Gamma[Q,D_i]}\in\F(H)$ for every $i$.
Equation~\eqref{eq:fundamental} gives \eqref{eq:intrinsic-analytic-criterion} on $\M\Omega$, hence on $H$ by boundedness and density.
Theorem~\ref{thm:intrinsic-analytic}\textup{(\ref{item:intrinsic-analytic-existence-criterion})} gives \textup{(1)}.

The rank independence is Lemma~\ref{lem:fraction-ranks}.
\end{proof}

\begin{proposition}\label{prop:fraction-commutator-rank}
Let $X$ and $(D_1,\ldots,D_d)$ satisfy the hypotheses of Theorem~\ref{thm:analytic}.
Let $(u,A,v)$ be a minimal linear representation of a nonzero $r\in\freefield{x_1,\ldots,x_d}$ as in Definition~\ref{def:linear-representation}, with $A=A_0+\sum_{i=1}^d A_ix_i$ of size $n$.
Then for every fraction pair $Q$ for $r(X)$,
\[
 \rank\overline{\Gamma[Q,D_i]}=\rank_{\C}A_i
 =\text{the tensor rank of }\partial_i r
 \qquad(1\leq i\leq d).
\]
Here, the tensor rank of a tensor $x$ is the least number $p$ such that $x=\sum_{\nu=1}^p y_\nu\otimes z_\nu$ for some $y_\nu$ and $z_\nu$.
\end{proposition}
\begin{proof}
Write $A=(A_{jk})_{j,k=1}^n$.
For $0\leq i\leq d$, let $A_{jk}^{(i)}$ be the $(j,k)$-entry of $A_i$.
Then
\[
 A_{jk}=A_{jk}^{(0)}+\sum_{i=1}^d A_{jk}^{(i)}x_i\in\C\langle x_1,\ldots,x_d\rangle.
\]

Set $g:=uA^{-1}$ and $h:=A^{-1}v$, with entries $g_j$ and $h_k$, respectively.
Applying $\partial_i$ entrywise to $Ah=v$ gives
\[
 0=(\partial_i A)h+A\partial_i h.
\]
Since $r=uA^{-1}v=uh$ and $\partial_i A_{jk}=A_{jk}^{(i)}(1\otimes1)$, we obtain
\[
 \partial_i r=u\partial_i h=-g(\partial_i A)h
 =-\sum_{j,k=1}^n A_{jk}^{(i)}\,g_j\otimes h_k^{\mathrm{op}}.
\]
By Theorem~\ref{thm:minimal-linear-representation}, both $(g_j)_{j=1}^n$ and $(h_k)_{k=1}^n$ are linearly independent over $\C$.
Choose a decomposition $\partial_i r=\sum_{\nu=1}^p y_\nu\otimes z_\nu^{\mathrm{op}}$ with the least number $p$ of simple tensors $y_\nu\otimes z_\nu^{\mathrm{op}}$.
By linear independence, choose linear functionals $\varphi_j$ on $g_j$ and $\widetilde{\varphi}_k$ on $h_k^{\mathrm{op}}$ such that $\varphi_j(g_{j'})=\delta_{jj'}$ and $\widetilde{\varphi}_k(h_{k'}^{\mathrm{op}})=\delta_{kk'}$.
Set $U_{j\nu}:=\varphi_j(y_\nu)$ and $V_{\nu k}:=\widetilde{\varphi}_k(z_\nu^{\mathrm{op}})$.
Applying $\varphi_j\otimes\widetilde{\varphi}_k$ to the two expressions for $\partial_i r$ gives $-A_{jk}^{(i)}=\sum_{\nu=1}^p U_{j\nu}V_{\nu k}$.
Thus $-A_i=UV$ and $\rank_{\C}A_i\leq p$.
Conversely, a rank factorization $-A_i=UV$ with $m:=\rank_{\C}A_i$ columns in $U$ and $m$ rows in $V$ gives
\[
 \partial_i r
 =\sum_{\nu=1}^m
 \left(\sum_{j=1}^n U_{j\nu}g_j\right)
 \otimes
 \left(\sum_{k=1}^n V_{\nu k}h_k^{\mathrm{op}}\right).
\]
Thus $p\leq m$, and the tensor rank of $\partial_i r$ is $m$.

Put $\M:=W^*(X)$.
Remark~\ref{rem:dual-maximality} gives $\Delta(X)=d$, so evaluation at $X$ is injective on the free field by Lemma~\ref{lem:free-field}.
Set $f:=r(X)$.
Apply Lemma~\ref{lem:clearing} to $f,g_1(X),\ldots,g_n(X)$ on the left and $f,h_1(X),\ldots,h_n(X)$ on the right.
We obtain regular $s,t\in\M$ such that
\[
 sf,\;ft,\;sg_j(X),\;h_k(X)t\in\M
 \qquad(1\leq j,k\leq n).
\]
Set $a:=sf$ and $b:=ft$.
Then $f=s^{-1}a=bt^{-1}$, so $Q:=(s,a,b,t)$ is a fraction pair for $f$.
Evaluating the formula for $\partial_i r$ in both tensor factors and multiplying by $s\otimes t^{\mathrm{op}}$ gives
\[
 \ell_i:=(s\otimes t^{\mathrm{op}})\partial_i f
 =-\sum_{j,k=1}^n A_{jk}^{(i)}\,sg_j(X)\otimes(h_k(X)t)^{\mathrm{op}}.
\]
By equations~\eqref{eq:psi-left-linearity} and~\eqref{eq:tensor-criterion}, these tensors satisfy
\[
 \sum_{i=1}^d\Psi_i(\ell_i)
 =(s\otimes t^{\mathrm{op}})(f\otimes1-1\otimes f^{\mathrm{op}})
 =a\otimes t^{\mathrm{op}}-s\otimes b^{\mathrm{op}}.
\]
Since all the $\ell_i$ lie in $\M\otimes_{\C}\M^{\mathrm{op}}$, Lemmas~\ref{lem:tensor-certificate} and~\ref{lem:finite-rank-identification} give $\overline{\Gamma[Q,D_i]}=-\Theta(\ell_i)$.

Define the maps
\[
 \begin{aligned}
 B&:\C^n\longrightarrow H, & Bc&:=\sum_{j=1}^n c_j\,sg_j(X)\Omega,\\
 C&:H\longrightarrow\C^n, & C\xi&:=\bigl(\langle J(h_k(X)t\Omega),\xi\rangle\bigr)_{k=1}^n .
 \end{aligned}
\]
Evaluation at $X$ and multiplication by the regular denominators preserve the linear independence of the two families.
Thus $B$ is injective, and $C$ is surjective because its adjoint $C^*$ sends the standard basis of $\C^n$ to the linearly independent vectors $J(h_k(X)t\Omega)$.
For every $\xi\in H$, equation~\eqref{eq:theta} gives
\begin{align*}
 \Theta(\ell_i)\xi
 &=-\sum_{j,k=1}^n A_{jk}^{(i)}\,
   \langle J(h_k(X)t\Omega),\xi\rangle\,sg_j(X)\Omega\\
 &=-\sum_{j=1}^n
   \left(\sum_{k=1}^n A_{jk}^{(i)}(C\xi)_k\right)sg_j(X)\Omega\\
 &=-BA_iC\xi.
\end{align*}
It follows that $\overline{\Gamma[Q,D_i]}=BA_iC$ has rank $\rank_{\C}A_i$, since $B$ is injective and $C$ is surjective.
Lemma~\ref{lem:fraction-ranks} gives the same rank for every fraction pair for $f$.
\end{proof}

When \(d=1\), let \(D_1\) be a dual operator for a normal operator \(X\) as in Proposition~\ref{prop:fraction-commutator-rank}.
Write $A(x)=A_0+xA_1$ for a minimal linear representation $(u,A,v)$ of $r$ used in Proposition~\ref{prop:fraction-commutator-rank}.
The determinant $\det A(x)$ is a nonzero polynomial because $A(x)$ is invertible over $\C(x)$.
Choose $\lambda\in\C$ with $A(\lambda)$ invertible.
Since $\det A(\lambda)\ne 0$, the entries of $A(x)^{-1}=\operatorname{adj}(A(x))/\det A(x)$ are analytic near $\lambda$, and so is $r(x)=uA(x)^{-1}v$.
The Taylor expansion of $r$ at $\lambda$ is
\[
 r(\lambda+z)=\sum_{n\geq 0}c_nz^n,
 \qquad c_n=\frac{r^{(n)}(\lambda)}{n!}.
\]
The free difference quotient is the divided difference, which has the expansion
\[
 \frac{r(\lambda+z)-r(\lambda+w)}{z-w}
 =\sum_{j,k\geq 0} c_{j+k+1}z^jw^k.
\]
Hence the tensor rank of $\partial_1 r$ equals the algebraic rank of $(c_{j+k+1})_{j,k\geq 0}$, the Hankel matrix of the Laurent series $r(\lambda+z^{-1})-r(\lambda)$.
Recall that Kronecker's theorem \cite[Articles~VII and~IX]{Kro81} characterizes rationality of a Laurent series vanishing at infinity by finite rankness of its Hankel matrix.
Thus Proposition~\ref{prop:fraction-commutator-rank} can be understood as a noncommutative analogue of Kronecker's theorem, relating the rank of the commutator to a measure of the rationality complexity of the underlying function.

Recall that the dimension of a minimal linear representation of a rational function $r$ is a purely algebraic quantity independent of the choice of realization by any embedding $r\mapsto r(X)$ for a tuple $X$ of operators.
The following corollary says we can give an upper bound for this dimension in terms of the ranks of the commutators $\overline{\Gamma[Q,D_i]}$.

\begin{corollary}\label{cor:minimal-representation-dimension}
Under the hypotheses of Proposition~\ref{prop:fraction-commutator-rank}, the dimension $n$ of a minimal linear representation of a rational function $r$ satisfies
\[
 n\leq 1+\sum_{i=1}^d\rank\overline{\Gamma[Q,D_i]}
\]
for every fraction pair $Q$ for $r(X)$.
\end{corollary}

\begin{proof}
Write $(u,A,v)$ for a minimal linear representation with $A=A_0+\sum_{i=1}^d A_ix_i$, and put $h:=A^{-1}v$.
If a complex row vector $w$ annihilates $\C v+\sum_{i=1}^d\Image A_i$, then $wv=0$ and $wA_i=0$ for every $1\leq i\leq d$, so
\[
 wA_0h=wAh=wv=0.
\]
By Theorem~\ref{thm:minimal-linear-representation}, the entries of $h$ are linearly independent over $\C$, so $wA_0=0$.
Thus $wA=0$ and $w=0$ by invertibility of $A$.
It follows that
\[
 \C^n=\C v+\sum_{i=1}^d\Image A_i,
 \qquad
 n\leq 1+\sum_{i=1}^d\rank_{\C}A_i.
\]
Proposition~\ref{prop:fraction-commutator-rank} now gives the asserted bound.
\end{proof}

\section{Examples}
\label{sec:examples}

\subsection{Scope of the rationality criteria}

The first two examples show that Theorems~\ref{thm:delta}, \ref{thm:intrinsic-analytic}, and \ref{thm:analytic} apply to successively smaller classes of generating tuples.

Our first example shows that maximal $\Delta$ does not imply the existence of a liberation dual system, even for selfadjoint generators whose scalar shifts are regular.

\begin{example}\label{ex:maximal-delta-without-liberation}
Let $(S_1,S_2)$ be a standard free semicircular pair in $(\mathcal N,\tau_{\mathcal N})$, with $\mathcal N=W^*(S_1,S_2)$.
In $\M:=\mathcal N\oplus\mathcal N$ with the trace $\tau:=\tfrac12(\tau_{\mathcal N}\oplus\tau_{\mathcal N})$, set
\[
 X_j:=S_j\oplus(S_j+5),\qquad j=1,2.
\]
We write $C^*(x_1,\ldots,x_m)$ for the $C^*$-subalgebra of $\M$ generated by $1,x_1,\ldots,x_m$.
The spectra of $S_j$ and $S_j+5$ lie in $[-2,2]$ and $[3,7]$, respectively.
Choose a continuous function $f:\mathbb R\to\mathbb R$ that equals $1$ on the first interval and $0$ on the second.
For each $j$, continuous functional calculus gives $f(X_j)\in C^*(X_j)$ and
\[
 f(X_j)=f(S_j)\oplus f(S_j+5)=1\oplus0=:e.
\]
Compressing by $e$ and $1-e$ recovers the generating pair on each summand, hence $\M=W^*(X_1,X_2)$.
The law of $X_j$ is the average of the law of $S_j$ and $S_j+5$, so it has no atoms.
Thus every scalar shift of $X_j$ is regular.

The free semicircular pair has the bounded dual system of right creation operators \cite[5.13]{Voi98}, so $\Delta(S_1,S_2)=2$ by \cite[Proposition~6.2]{MSY23}.
Let $A$ be a full square linear matrix over $\C\langle x_1,x_2\rangle$.
The translation $x_j\mapsto x_j+5$ is an automorphism of the free algebra, so $A(x_1+5,x_2+5)$ is also full.
By \cite[Theorem~3.1]{MSY23}, both $A(S_1,S_2)$ and $A(S_1+5,S_2+5)$ have zero kernel.
Their direct sum is $A(X_1,X_2)$ and has zero kernel as well.
Applying the same theorem to $(X_1,X_2)$ gives $\Delta(X_1,X_2)=2$.

Suppose that $(L_1,L_2)$ is a liberation dual system as in Definition~\ref{def:liberation-dual-system}.
On $\M\Omega$, its commutator relations give
\[
 [X_1,L_1-P_\Omega]=0,\qquad [X_2,L_1]=0.
\]
Choose real polynomials $p_n$ converging uniformly to $1$ on $[-2,2]$ and to $0$ on $[3,7]$, so that $p_n(X_j)\rightarrow e$ in operator norm for both $j$.
For either $(Z,B)=(X_1,L_1-P_\Omega)$ or $(Z,B)=(X_2,L_1)$, commutation with $Z$ gives $p_n(Z)B\Omega=Bp_n(Z)\Omega$.
For $\xi\in\M\Omega$, pairing the latter vector with $\xi$ gives
\[
 \langle \xi,p_n(Z)B\Omega\rangle
 =\langle B^*\xi,p_n(Z)\Omega\rangle
 \longrightarrow\langle B^*\xi,e\Omega\rangle
 =\langle \xi,Be\Omega\rangle .
\]
Since $p_n(Z)B\Omega\rightarrow eB\Omega$ and $\M\Omega$ is dense, $eB\Omega=Be\Omega$.
The two choices of $B$ give
\[
 [e,L_1-P_\Omega]\Omega=0,\qquad [e,L_1]\Omega=0.
\]
Subtracting the first equality from the second yields $[e,P_\Omega]\Omega=0$.
This contradicts $[e,P_\Omega]\Omega=(e-\tfrac12)\Omega\ne0$.
\end{example}

The following example shows that the existence of a liberation dual system does not imply the existence of a dual system, even for normal generators with nonatomic spectral distributions.

\begin{example}\label{ex:liberation-without-duals}
Let
\[
 d\mu(t):=\frac{dt}{\pi\sqrt{1-t^2}}\quad(-1<t<1),\qquad
 \M:=L^\infty((-1,1),\mu),\qquad \tau(f):=\int f\,d\mu,
\]
and let $X(x):=x$.
Then $X=X^*$ generates $\M$ and has nonatomic spectral distribution, so every $X-\lambda$, $\lambda\in\C$, is regular.
Taking $L:=P_\Omega$ gives a normalized liberation dual system, since $[X,L]=[X,P_\Omega]$ and $L$ is bounded on $H$.

Suppose that $D$ is a dual system for $X$ in the sense of Definition~\ref{def:dual-system}.
The commutator product rule gives on $\M\Omega$
\[
 [X^{2n+1},D]=\sum_{k=0}^{2n}X^kP_\Omega X^{2n-k}.
\]
Write $m_r:=\tau(X^{2r})=4^{-r}\binom{2r}{r}$ for $r\geq0$.
Pairing the commutator identity with $\Omega$ gives
\[
 \begin{aligned}
 \langle\Omega,[X^{2n+1},D]\Omega\rangle
 &=\sum_{k=0}^{2n}\tau(X^k)\tau(X^{2n-k})\\
 &=\sum_{r=0}^{n}m_r m_{n-r},
 \end{aligned}
\]
where the last equality is due to $\tau(X^{2r+1})=0$.
We know that $\sum_{r=0}^\infty m_r z^r=(1-z)^{-1/2}$ for $|z|<1$, so the Cauchy product of this series with itself gives
\[
 \sum_{n=0}^\infty\Bigl(\sum_{r=0}^n m_r m_{n-r}\Bigr)z^n
 =\Bigl(\sum_{r=0}^\infty m_r z^r\Bigr)^2
 =(1-z)^{-1}
 =\sum_{n=0}^\infty z^n.
\]
Comparing coefficients of $z^n$ shows that $\langle\Omega,[X^{2n+1},D]\Omega\rangle=1$.

Since $X$ is the coordinate function, we have
\[
 \|X^{2n+1}\|_2^2=\int_{-1}^1 |t|^{4n+2}\,d\mu(t).
\]
For every $t\in(-1,1)$, $|t|^{4n+2}\to0$ and $|t|^{4n+2}\leq1$.
Since $\mu$ is a probability measure, dominated convergence implies $\|X^{2n+1}\|_2^2\to0$ and hence $\|X^{2n+1}\|_2\to0$.
Finally, we have
\[
 \langle\Omega,[X^{2n+1},D]\Omega\rangle
 =\langle X^{2n+1}\Omega,D\Omega\rangle
  -\langle D^*\Omega,X^{2n+1}\Omega\rangle.
\]
By Cauchy--Schwarz,
\[
 1=\bigl|\langle\Omega,[X^{2n+1},D]\Omega\rangle\bigr|
 \leq\bigl(\|D\Omega\|_2+\|D^*\Omega\|_2\bigr)\|X^{2n+1}\|_2
 \longrightarrow0,
\]
a contradiction.
Thus this tuple satisfies the hypotheses of Theorem~\ref{thm:intrinsic-analytic} but not those of Theorem~\ref{thm:analytic}.
\end{example}

Replacing the dual system in Theorem~\ref{thm:analytic} by a normalized liberation dual system can make the finite-rank condition hold for nonrational elements.
In the above example, $L=P_\Omega$ for an arcsine generator $X$, so for every $u\in L^0(\M)$ and every fraction pair $Q=(s,a,b,t)$ for $u$, both terms in
\[
 \Gamma[Q,L]=sP_\Omega b-aP_\Omega t
\]
have rank at most one, and hence $\rank\Gamma[Q,L]\leq 2$.
Since $\mu$ is nonatomic, $\Div_{L^0(\M)}\C[X]$ is the field of rational functions of $X$.
The projection $p:=\mathbf 1_{(0,1)}(X)$ is nonrational, but we cannot use the finite-rank condition of $\Gamma[Q,L]$ to detect this.
For this projection, equation~\eqref{eq:intrinsic-analytic-coordinates} takes the form
\[
 [\rho(X),T]=\Gamma[Q,L].
\]
By Theorem~\ref{thm:intrinsic-analytic}, it has no solution $T\in\F(H)$ for any fraction pair $Q$ for $p$.

For a fixed bounded rational element, the finite-rank commutator equation~\eqref{eq:certificate} in Theorem~\ref{thm:delta} can have a solution for one fraction pair and no solution for another.

\begin{example}\label{ex:choice-fraction-pair}
Let $h$ be selfadjoint with the uniform distribution on $(0,1)$, let $S$ be a standard semicircular element free from $h$, and put $\M:=W^*(h,S)$ and $Y:=hSh$.
The operator $h$ is regular, so $h^{-1}\in L^0(\M)$, while the integral $\int_0^1 t^{-2}\,dt=\infty$ gives $h^{-1}\notin H:=L^2(\M)$.
We note that $S=h^{-1}Yh^{-1}$ is rational in $(h,Y)$.

The distributions of both $h$ and $S$ has no atoms, and hence $\Delta(h,S)=2$ by \cite[Section~6.1]{MSY23}.
Let $T_h,T_Y$ be finite-rank operators.
Applying the commutator product rule to $\rho(Y)=\rho(h)\rho(S)\rho(h)$ gives
\[
 \begin{aligned}
 [\rho(Y),T_Y]
 &=[\rho(h),\rho(S)\rho(h)T_Y]\\
 &\quad+[\rho(S),\rho(h)T_Y\rho(h)]\\
 &\quad+[\rho(h),T_Y\rho(h)\rho(S)].
 \end{aligned}
\]
Adding $[\rho(h),T_h]$ and using $\rho(S)\rho(h)=\rho(hS)$ and $\rho(h)\rho(S)=\rho(Sh)$, we obtain
\begin{equation}\label{eq:choice-fraction-pair-commutators}
 \begin{aligned}
 [\rho(h),T_h]+[\rho(Y),T_Y]
 &= [\rho(h),T_h+\rho(hS)T_Y+T_Y\rho(Sh)]\\
 &\quad+[\rho(S),\rho(h)T_Y\rho(h)].
 \end{aligned}
\end{equation}
If the left-hand side of equation~\eqref{eq:choice-fraction-pair-commutators} vanishes, the maximality $\Delta(h,S)=2$ gives
\[
\rho(h)T_Y\rho(h)=0 \quad\text{and}\quad T_h+\rho(hS)T_Y+T_Y\rho(Sh)=0.
\]
Since $\rho(h)$ is injective and has dense range, we get $T_Y=0$ and then $T_h=0$.
Thus $\Delta(h,Y)=2$ by Remark~\ref{rem:cycle-characterization} and thus Theorem~\ref{thm:delta} applies to $(h,Y)$.

For the fraction pair $(1,S,S,1)$, the right hand side of equation~\eqref{eq:certificate} is
\begin{equation}\label{eq:choice-fraction-pair-commutator-S}
[P_\Omega,S]=-[\rho(S),P_\Omega].
\end{equation}
If $(T_h,T_Y)$ were a finite-rank solution of \eqref{eq:certificate}, using equation~\eqref{eq:choice-fraction-pair-commutators} and \eqref{eq:choice-fraction-pair-commutator-S} we obtain
\[
 [\rho(h),T_h+\rho(hS)T_Y+T_Y\rho(Sh)]+[\rho(S),\rho(h)T_Y\rho(h)+P_\Omega]=0.
\]
The maximality $\Delta(h,S)=2$ forces $\rho(h)T_Y\rho(h)=-P_\Omega$.
Applying this identity to $\Omega$ yields $\rho(h)T_Y(h\Omega)=-\Omega$, so $T_Y(h\Omega)=-h^{-1}$ and this forces $h^{-1}\in H$.
This contradicts $h^{-1}\notin H$.

For $Q=(h,hS,Sh,h)$, take $T_Y=-P_\Omega$ and $T_h=\rho(hS)P_\Omega+P_\Omega\rho(Sh)$.
The first commutator on the right of equation~\eqref{eq:choice-fraction-pair-commutators} vanishes, and
\[
 [\rho(h),T_h]+[\rho(Y),T_Y]
 =-[\rho(S),hP_\Omega h]
 =hP_\Omega S h-hS P_\Omega h
 =\Gamma[Q,P_\Omega].
\]
\end{example}

The last example in this subsection illustrates the necessity of the normality assumption in Theorem~\ref{thm:analytic}.

\begin{example}\label{ex:normality}
Let $l_j,r_j$, $j=1,2$, be the left and right creation operators on a full Fock space with the vacuum vector $\Omega$.
Their adjoints $l_j^*,r_j^*$ are the left and right annihilation operators, respectively.
The operators $S_j:=l_j+l_j^*$ form a standard free semicircular pair with respect to the vacuum trace $\tau(x):=\langle\Omega,x\Omega\rangle$ on $W^*(S_1,S_2)$ (see \cite[5.13]{Voi98}).
Put
\[
 c:=\frac{S_1+\mathrm iS_2}{\sqrt2},\qquad
 \M:=W^*(c)=W^*(S_1,S_2).
\]
Let us consider $c^*$ as an element of $L^0(\M)$.
Since every element of $\Div_{L^0(\M)}\C[c]$ commutes with $c$ and $[c,c^*]=\mathrm i[S_2,S_1]\ne0$, we have $c^*\notin\Div_{L^0(\M)}\C[c]$.

Define
\[
 D:=-\frac{r_1^*-\mathrm i r_2^*}{\sqrt2}\in\B(H).
\]
The identities $[r_i^*,S_j]=\delta_{ij}P_\Omega$ give $[c,D]=P_\Omega$ and $[D,c^*]=0$.
Thus, for every fraction pair $Q=(s,a,b,t)$ for $c^*$,
\[
 \Gamma[Q,D]=s[D,c^*]t=0\in\F(H).
\]
Thus the implication from condition~\textup{(\ref{item:analytic-fraction-commutators})} of Theorem~\ref{thm:analytic} to rationality fails when normality is dropped.

The bounded dual $D$ gives $\Delta(c)=1$ by \cite[Proposition~6.2]{MSY23}.
The free-field embedding in Lemma~\ref{lem:free-field} then shows that $c-\lambda$ is invertible in $L^0(\M)$ for every $\lambda\in\C$.
The single operator $L:=P_\Omega$ is a normalized liberation dual system for $c$, so Theorem~\ref{thm:intrinsic-analytic} applies.
For \(u=c^*\), the commutator equation in assertion~\eqref{item:intrinsic-analytic-criterion} has no finite-rank solution for any fraction pair.
Fix a fraction pair $Q=(s,a,b,t)$ for $c^*$, and suppose that $T\in\F(H)$ satisfies equation~\eqref{eq:intrinsic-analytic-criterion}.
Lemma~\ref{lem:finite-rank-identification} then gives $T=\overline{\Gamma[Q,D]}=0$, so
\[
 0=[\rho(c),T]=\Gamma[Q,P_\Omega]=s[P_\Omega,c^*]t.
\]
Since $\tau(c^*)=0$, we have
\[
 [P_\Omega,c^*]\Omega=-c^*\Omega\ne0,
\]
where $c^*\Omega$ is nonzero by faithfulness of the vacuum trace.
Regularity makes $s$ injective and $t$ have dense range, so $s[P_\Omega,c^*]t\ne0$, a contradiction.

\end{example}

\subsection{Division closures in bounded operators}\label{subsec:bounded-division-closures}

For a tuple $X$ of free Haar unitaries \cite[Theorem~1.1]{Lin00} or of free semicircular elements \cite[Theorem~4.5]{Miy22}, set $\M:=W^*(X)$ and $\A_X:=\C\langle X\rangle$.
Then
\[
 \Div_{\M}\A_X
 =\Div_{L^0(\M)}\A_X\cap\M.
\]
This equality can fail even for selfadjoint generators admitting a bounded dual system.
At moment, we do not have a general criterion for when this equality holds.

To give a concrete example that fails the equality, we first introduce the following proposition.
This proposition may be known to experts, but we include a proof for completeness.
An alternative argument can also be found based on the results by Pincus \cite[Theorem~2.1 and p.~226]{Pin68} or by Clancey and Putnam \cite[p.~452]{CP71}.

\begin{proposition}\label{prop:scalar-bounded-dual}
Let $X$ be a bounded selfadjoint element of a tracial $W^*$-probability space $(\M,\tau)$ and let $\mu$ be its distribution.
Then $X$ admits a bounded dual system $(D)$ on $L^2(\M,\tau)$ if and only if $d\mu(t)=w(t)\,dt$ for some $w\in L^\infty(\mathbb R)$.
In this case, the dual operator can be chosen such that
\[D^*=-D,\qquad\text{and}\qquad\|D\|\leq\pi\|w\|_\infty.\]
\end{proposition}

\begin{proof}
Identify $K:=L^2(W^*(X),\tau)$ with $L^2(\mathbb R,\mu)$, where $X$ acts by coordinate multiplication, and let $P_1$ be the orthogonal projection onto $\C1$.
For every bounded interval $I$ of positive length, the spectral projection $E_X(I)$ satisfies
\[
 \|P_1 E_X(I)\|^2=\|E_X(I)1\|_2^2=\mu(I).
\]
Let $\|\cdot\|_X$ denote Kato's $X$-smoothness norm \cite[(2.2), p.~536]{Kat68}.
The preceding identity gives
\[
 c:=\|P_1\|_X^2=\sup_{0<|I|<\infty}\frac{\mu(I)}{|I|},
\]
where the supremum is over intervals $I$.
Thus $P_1$ is $X$-smooth if and only if $c<\infty$.
If $c<\infty$, covering Borel sets by intervals gives $\mu\leq c\,dt$.
Thus $\mu$ has a density $w$ with $w\in L^\infty(\mathbb R)$.
Conversely, if $d\mu(t)=w(t)\,dt$ with $w\in L^\infty(\mathbb R)$, then $\mu(I)\leq\|w\|_\infty |I|$ for every such interval.
Thus $c=\|w\|_\infty$ whenever the density is essentially bounded.

Since $P_1^{1/2}=P_1$, Kato's commutator theorem \cite[Theorem~(6.2), p.~543]{Kat68} states that $P_1$ is $X$-smooth if and only if there exists a bounded selfadjoint operator $B$ with $\mathrm i[X,B]=P_1$.
For such a $B$, the operator $D:=\mathrm iB$ satisfies $[X,D]=P_1$.
Conversely, a bounded $D$ with $[X,D]=P_1$ satisfies $[X,D^*]=-P_1$, so $B:=(D-D^*)/(2\mathrm i)$ is a bounded selfadjoint solution of $\mathrm i[X,B]=P_1$.
The proof of Kato's theorem \cite[p.~545]{Kat68} gives a solution $B$ with $0\leq B\leq 2\pi c I_K$.
Hence $D:=\mathrm i(B-\pi c I_K)$ is the desired dual operator.
\end{proof}

\begin{example}\label{ex:bounded-division-closure}
In a tracial $W^*$-probability space, let $X$ be a positive operator with distribution $d\mu(t)=4t^3\,dt$ on $(0,1)$, and let $S$ be a standard semicircular element free from $X$.
Put
\[
 Y:=XSX,\qquad \M:=W^*(X,Y),\qquad \A:=\C\langle X,Y\rangle.
\]
The selfadjoint pair $(X,Y)$ admits a bounded dual system.
For this pair,
\[
 \Div_{\M}\A
 \subsetneq \Div_{L^0(\M)}\A\cap\M.
\]
The strict inclusion is guaranteed by $S=X^{-1}YX^{-1}\in\Div_{L^0(\M)}\A\cap\M$ and $S\notin\Div_{\M}\A$.

We prove the statement in the following steps.
Since the law of $X$ has no atom at zero and the trace is faithful, $X$ has zero kernel.
It implies that $X^{-1}\in L^0(\M)$ and thus
\[
 S=X^{-1}YX^{-1}\in\Div_{L^0(\M)}\A.
\]
So we see that $S\in L^0(\M)$ is bounded and affiliated with $\M$, it therefore has to be in $\M$ \cite[Proposition~7.2.3]{AP18}.

Next, we want to show that $S\notin C^*(X,Y)$.
For that purpose, put
\[
 \pi_\varepsilon:=1_{(0,\varepsilon]}(X),\qquad 0<\varepsilon<1.
\]
The distribution of $X$ gives
\[
 \tau(\pi_\varepsilon)=\mu((0,\varepsilon])
 =\int_0^\varepsilon 4t^3\,dt=\varepsilon^4>0.
\]
Define
\[
 \xi_\varepsilon:=\frac{\pi_\varepsilon\Omega}{\sqrt{\tau(\pi_\varepsilon)}}.
\]
Since $\|\pi_\varepsilon\Omega\|_2^2=\tau(\pi_\varepsilon)$, we have $\|\xi_\varepsilon\|_2=1$.
Moreover,
\[
 \|X\xi_\varepsilon\|_2\leq\varepsilon,\qquad
 \|Y\xi_\varepsilon\|_2\leq 2\varepsilon,
\]
since $\|X\|=1$ and $\|S\|=2$.
For a nonconstant word $Z_1\cdots Z_m$ with $m\geq1$ and $Z_j\in\{X,Y\}$, these estimates give
\[
 \begin{aligned}
 \|Z_1\cdots Z_m\xi_\varepsilon\|_2
 &\leq\|Z_1\cdots Z_{m-1}\|\,\|Z_m\xi_\varepsilon\|_2\\
 &\leq2\varepsilon\|Z_1\cdots Z_{m-1}\|
 \longrightarrow0\qquad(\varepsilon\rightarrow0).
 \end{aligned}
\]
Thus for a polynomial $p$ with constant term $c=p(0,0)$ we have $\|(p(X,Y)-c)\xi_\varepsilon\|_2\to0$.
Put $B:=(S-c)^*(S-c)$, we see that
\[
 \begin{aligned}
 \|(S-c)\xi_\varepsilon\|_2^2
 &=\frac{\|(S-c)\pi_\varepsilon\Omega\|_2^2}{\tau(\pi_\varepsilon)}
 =\frac{\tau(\pi_\varepsilon B\pi_\varepsilon)}{\tau(\pi_\varepsilon)}\\
 &=\frac{\tau(B\pi_\varepsilon^2)}{\tau(\pi_\varepsilon)}
 =\frac{\tau(B\pi_\varepsilon)}{\tau(\pi_\varepsilon)}\\
 &=\frac{\tau(B)\tau(\pi_\varepsilon)}{\tau(\pi_\varepsilon)}
 =\tau(B)\\
 &=\tau(S^2)-(c+\overline{c})\tau(S)+|c|^2
 =1+|c|^2,
 \end{aligned}
\]
where the last equality uses $\tau(S)=0$ and $\tau(S^2)=1$.
The reverse triangle inequality gives
\[
 \bigl|\|(S-p(X,Y))\xi_\varepsilon\|_2-\sqrt{1+|c|^2}\bigr|
 \leq\|(p(X,Y)-c)\xi_\varepsilon\|_2
 \xrightarrow{\varepsilon\rightarrow0}0.
\]
Since each $\xi_\varepsilon$ is a unit vector,
\[
 \|S-p(X,Y)\|\geq\lim_{\varepsilon\rightarrow0}\|(S-p(X,Y))\xi_\varepsilon\|_2=\sqrt{1+|c|^2}\geq1.
\]
Since $X,Y$ are selfadjoint, $\A$ is norm dense in $C^*(X,Y)$, so $\dist(S,C^*(X,Y))\geq1$.
The unital $C^*$-subalgebra $C^*(X,Y)$ is closed under inversion in $\M$ \cite[Corollary~II.1.6.7]{Bla06}, so $\Div_{\M}\A=\Div_{C^*(X,Y)}\A$.
Since $S\notin C^*(X,Y)$, we have $S\notin\Div_{\M}\A$.

It remains to construct a bounded dual system for $(X,Y)$.
Since $S\in\M=W^*(X,Y)$ and $Y=XSX$, the free pair $(X,S)$ generates $\M$.
Let $H_0:=L^2((0,1),\mu)$ and let $x\in H_0$ be the coordinate function $x(t)=t$.
We consider $X,Y,S$ as operators acting on
\[
 H=L^2(\M,\tau)=\bigoplus_{n\geq0}H_n,\qquad H_n:=H_0^{\bar\otimes(n+1)}.
\]
The vacuum vector $\Omega$ is the constant function $1$ in $H_0$, and $X$ acts by coordinate multiplication $x$ on the first factor in each $H_n$.
Let $l$ be the isometry determined by
\[
 l(\xi_0\otimes\cdots\otimes\xi_n)
 =1\otimes\xi_0\otimes\cdots\otimes\xi_n.
\]
Its adjoint $l^*$ vanishes on $H_0$ and is given by
\[
 l^*(\xi_0\otimes\cdots\otimes\xi_n)
 =\langle 1,\xi_0\rangle\,\xi_1\otimes\cdots\otimes\xi_n.
\]
Then we have $S=l+l^*$ on $H$.
Let $L:=XlX$ and then we have
\[
 L(\xi_0\otimes\cdots\otimes\xi_n)
 =x\otimes X\xi_0\otimes\xi_1\otimes\cdots\otimes\xi_n.
\]
with its adjoint $L^*$ given by
\[
 L^*(\xi_0\otimes\cdots\otimes\xi_n)
 =\langle x,\xi_0\rangle X\xi_1\otimes\xi_2\otimes\cdots\otimes\xi_n.
\]
Thus
\[
Y=XSX=X(l+l^*)X=XlX+Xl^*X=L+L^*.
\]

Define $K$ from $H_0\bar\otimes H_0$ to $H_0$ by the kernel
\[
 k(s,t):=\frac{1_{\{t\leq s\}}}{st},\qquad
 (K\zeta)(s):=\int_0^1 k(s,t)\zeta(s,t)\,d\mu(t).
\]
For $s\in(0,1)$,
\[
 \int_0^1 |k(s,t)|^2\,d\mu(t)
 =\frac{4}{s^2}\int_0^s t\,dt=2.
\]
Thanks to Cauchy--Schwarz equality, we obtain
\[
 \|K\zeta\|_2^2
 \leq 2\int_0^1\int_0^1|\zeta(s,t)|^2\,d\mu(t)\,d\mu(s)
 =2\|\zeta\|_2^2.
\]
Let $\tilde{R}$ be a linear operator on $H$ that vanishes on $H_0$ and
\[
\tilde{R}|_{H_n}:=I_{H_0}^{\otimes(n-1)}\otimes K\quad\text{for } n\geq1.
\]
Then we have $\|\tilde{R}\|\leq\|K\|\leq\sqrt{2}$.
Note that the formula for $K$ gives $K(X\otimes I_{H_0})=XK$, so we see that $\tilde{R}$ commutes with $X$ on $H$ by a simple verification on each $H_n$.

Next, we calculate $\tilde{R}L$ and $L\tilde{R}$ on $H_0$ and and show that $\tilde{R}L=L\tilde{R}$ on $H_n$ for $n\geq1$.
For $f\in H_0$, $L\tilde{R}f=0$ and
\[
 (\tilde{R}Lf)(s)=\bigl(K(x\otimes Xf)\bigr)(s)
 =s\int_0^1 k(s,t)t f(t)\,d\mu(t)
 =\int_0^s f(t)\,d\mu(t).
\]
Write $V$ for this integral operator on $H_0$.
On $H_1$, the identity $K(X\otimes I_{H_0})=XK$ gives
\[
 \tilde{R}L(\xi_0\otimes\xi_1)
 =x\otimes K(X\xi_0\otimes\xi_1)
 =x\otimes X K(\xi_0\otimes\xi_1)
 =L\tilde{R}(\xi_0\otimes\xi_1).
\]
On $H_n$ for $n\geq2$, $L$ prepends $x$ and multiplies the first factor by $X$, while $\tilde{R}$ acts on the last two factors.
Thus $\tilde{R}L=L\tilde{R}$ on $H_n$ for $n\geq2$.

We next verify that $\tilde{R}L^*=L^*\tilde{R}$ on $H_n$ for $n\geq0$.
Both $\tilde{R}L^*$ and $L^*\tilde{R}$ vanish on $H_0$ and $H_1$ by definition.
On $H_2$ we have
\[
 \begin{aligned}
 \tilde{R}L^*(\xi_0\otimes\xi_1\otimes\xi_2)
 &=\langle x,\xi_0\rangle K(X\xi_1\otimes\xi_2)\\
 &=\langle x,\xi_0\rangle X K(\xi_1\otimes\xi_2)
 =L^*\tilde{R}(\xi_0\otimes\xi_1\otimes\xi_2).
 \end{aligned}
\]
On $H_n$ for $n\geq3$, $L^*$ acts on the first two factors and $\tilde{R}$ on the last two, so the two compositions agree there as well.

As $Y=L+L^*$, the above computations give
\begin{equation}\label{eq:sandwich-last-two-factor-commutator}
 [\tilde{R},Y]=V\oplus0\qquad\text{on } H_0\oplus\bigoplus_{n\geq1}H_n.
\end{equation}
Set $D_2:=\tilde{R}^*-\tilde{R}$.
Since $X=X^*$, taking adjoints in $\tilde{R}X=X\tilde{R}$ gives $X\tilde{R}^*=\tilde{R}^*X$.
Thus both terms in $[X,D_2]=[X,\tilde{R}^*]-[X,\tilde{R}]$ vanish.
Since $Y=Y^*$,
\[
 [\tilde{R},Y]^*=(\tilde{R}Y-Y\tilde{R})^*=Y\tilde{R}^*-\tilde{R}^*Y=[Y,\tilde{R}^*].
\]
Using equation~\eqref{eq:sandwich-last-two-factor-commutator}, we obtain
\[
 [Y,D_2]=[Y,\tilde{R}^*]-[Y,\tilde{R}]=[\tilde{R},Y]^*+[\tilde{R},Y]=(V^*+V)\oplus0.
\]
The kernel of $V$ is $1_{\{t\leq s\}}$, so $(V^*f)(s)=\int_s^1 f(t)\,d\mu(t)$.
Since $\mu$ has no atoms,
\[
 ((V^*+V)f)(s)=\int_0^1 f(t)\,d\mu(t)=\langle 1,f\rangle.
\]
So we see that $V^*+V=P_{\Omega}|_{H_0}$.
Consequently, we obtain
\begin{equation}\label{eq:sandwich-second-dual}
 [X,D_2]=0,\qquad [Y,D_2]=P_\Omega.
\end{equation}

The density of $\mu$ has essential supremum $4$, so Proposition~\ref{prop:scalar-bounded-dual} gives a operator $D_X$ on $H_0$ with $[X,D_X]=P_\Omega|_{H_0}$ and $\|D_X\|\leq4\pi$.
Define the bounded skew-adjoint operator $\tilde{D}_X$ on $H$ by $\tilde{D}_X|_{H_n}:=I_{H_0}^{\otimes n}\otimes D_X$.
The commutator satisfies $[X,\tilde{D}_X]|_{H_0}=P_\Omega|_{H_0}$, and on every positive grade the actions of $X$ and $\tilde{D}_X$ are on distinct factors.
It follows that $[X,\tilde{D}_X]=P_\Omega$.

Set $Z:=YX^{-1}$ in $L^0(\M)$.
The identity $Y=XSX$ gives $Z=XS\in\M$.
Hence first we observe that for $f\in H_0$,
\[
 [L,\tilde{D}_X]f=x\otimes(X D_X-D_X X)f
 =x\otimes (P_\Omega|_{H_0})f
 =ZP_\Omega f.
\]
On $H_n$ for $n\geq1$, $L$ and $\tilde{D}_X$ commute because $L$ prepends x and multiplies the first factor by X, while $\tilde{D}_X$ acts on the last.
Hence $[L,\tilde{D}_X]=ZP_\Omega$.
Since $\tilde{D}_X^*=-\tilde{D}_X$ by Proposition~\ref{prop:scalar-bounded-dual}, taking adjoints gives $[L^*,\tilde{D}_X]=[L,\tilde{D}_X]^*$.
Therefore $[Y,\tilde{D}_X]=ZP_\Omega+P_\Omega Z^*$.

Finally, let us define the bounded operator
\[
 D_1:=\tilde{D}_X-\rho(Z)D_2-D_2\rho(Z^*),
\]
where $\rho(Z)$ denotes the right multiplication operator by $Z$.
Thus $[X,D_1]=[X,\tilde{D}_X]=P_\Omega$, thanks to equation~\eqref{eq:sandwich-second-dual}.
Since $\rho(Z)\Omega=Z\Omega$, we have $\rho(Z)P_\Omega=ZP_\Omega$.
For $b\in\M$, traciality gives
\[
 P_\Omega\rho(Z^*)b\Omega
 =\tau(bZ^*)\Omega
 =\tau(Z^*b)\Omega
 =P_\Omega Z^*b\Omega.
\]
Both sides are bounded operators, so $P_\Omega\rho(Z^*)=P_\Omega Z^*$ by density of $\M\Omega$.
Consequently
\[
 \begin{aligned}
 [Y,D_1]
 &=[Y,\tilde{D}_X]-\rho(Z)[Y,D_2]-[Y,D_2]\rho(Z^*)\\
 &=ZP_\Omega+P_\Omega Z^*-\rho(Z)P_\Omega-P_\Omega\rho(Z^*)\\
 &=0.
 \end{aligned}
\]
Thus $(D_1,D_2)$ is a bounded dual system for $(X,Y)$ in the sense of Definition~\ref{def:dual-system}.
\end{example}

\bibliographystyle{f0001-amsalpha}
\bibliography{references}
\end{document}